\documentclass[numbers,webpdf,imanum]{ima-authoring-template}%

\usepackage{amsmath}
\usepackage{amssymb}
\usepackage{graphicx} 
\usepackage{xcolor}
\usepackage{xspace}
\usepackage{multirow}
\usepackage{bbm}
\usepackage{bm}
\usepackage{Macros/mydef}
\usepackage{Macros/remarks}
\usepackage{algorithm}
\usepackage[numbers]{natbib}
\usepackage{graphicx}
\usepackage{multirow}
\usepackage{mathtools}
\usepackage{rotating}
\usepackage{pdflscape}
\usepackage[normalem]{ulem}

\graphicspath{{Fig/}}

\usepackage[mathlines]{lineno}

\theoremstyle{thmstyletwo}%
\newtheorem{theorem}{Theorem}
\newtheorem{proposition}[theorem]{Proposition}%
\newtheorem{lemma}[theorem]{Lemma}%

\numberwithin{equation}{section}

\newcommand{\reviewerOne}[1]{{\leavevmode\color{black}{#1}}}
\newcommand{\reviewerTwo}[1]{{\leavevmode\color{black}{#1}}}

\newcommand{\igor}[1]{\textcolor{black}{#1}}
\newcommand{\igorr}[1]{\textcolor{black}{#1}}
\newcommand{\igorrr}[1]{\textcolor{black}{#1}}

\newcommand{\surf}{\reviewerTwo{s}}
\newcommand{\edge}{\reviewerTwo{J}}

\begin{document}

\copyrightyear{2025, submitted}
\vol{}
\appnotes{Paper}
\firstpage{1}


\title{\reviewerTwo{Free‑surface Stokes problem: stability estimates and time‑step improvements}}

\author{Igor Tominec$^{*,1,2,4}$, Lukas Lundgren$^{1,2,3}$, André Löfgren$^{1,2}$, Josefin Ahlkrona$^{1,2,3}$
\vspace{0.2cm}
\address{\orgdiv{$^1$ Department of Mathematics}, \orgname{Stockholm University}, \orgaddress{\postcode{SE-106 91}, \state{Stockholm}, \country{Sweden}}} 
\vspace{-0.15cm}
\address{\orgdiv{$^2$ Bolin Centre for Climate Research}, \orgname{Stockholm University}, \orgaddress{\postcode{SE-106 91}, \state{Stockholm}, \country{Sweden}}}
\vspace{-0.15cm}
\address{\orgdiv{$^3$ Swedish e-Science Research Centre}, \orgaddress{\postcode{SE-100 44}, \state{Stockholm}, \country{Sweden}}}
\vspace{-0.15cm}
\address{\orgdiv{$^4$ Faculty of Mechanical Engineering, Laboratory for Fluid Dynamics and Thermodynamics}, \orgname{University of Ljubljana,} \orgaddress{\street{A\v{s}ker\v{c}eva cesta 6}, \postcode{1000}, \state{Ljubljana}, \country{Slovenia}}}
\vspace{0.3cm}
}


\authormark{Tominec, Lundgren, Löfgren, Ahlkrona}

\corresp[*]{Corresponding author: \href{email:igor.tominec@gmail.com}{igor.tominec@gmail.com}}

\received{12}{06}{2025}


\abstract{Accurate simulations of ice sheet dynamics, mantle convection, \igor{lava flow}, 
  and other highly viscous free-surface flows involve solving the coupled \igor{Stokes/free-surface equations}. 
  In this paper, we theoretically analyze the stability and conservation properties of the weak form of this system 
  for Newtonian fluids and non-Newtonian fluids, 
  at both the continuous and discrete levels. We perform the fully discrete stability analysis for 
  finite element methods used in space with explicit and implicit Euler time-stepping methods used in time. 
  Motivated by the theory, we propose a stabilization term designed for the explicit Euler discretization, which ensures unconditional time stability 
  \reviewerTwo{with respect to the time-step size 
  when free-surface height is Lipschitz continuous}, and permits conservation of the domain volume. 
  Numerical experiments validate and support our theoretical findings.}


\maketitle
\section{Introduction}
Free-surface flows of Newtonian and non-Newtonian fluids in the low Reynolds number regime are commonly modelled by the Stokes equations coupled with a free-surface evolution equation, 
resulting in a moving domain problem.  
Real-world applications include ice sheet dynamics \cite{greve_book}, mantle convection \cite{Kaus_fssa}, 
\igor{lava flow \cite{2020_hinton_lava}}, 
interactions between viscous drops \cite{Manga_visc_drops}, and others. 
These applications demand numerically robust simulations that remain stable across a range of mesh resolutions and time step sizes. 
However, standard discretizations of the coupled Stokes/free-surface problem may suffer from instabilities that 
constrain the simulations to small time step sizes, which then induce large computational costs. 
In this paper, we rigorously analyze the stability of the coupled Stokes/free-surface system and, motivated by our theoretical results, 
propose new numerical stabilization approaches. These enable stable and efficient simulations of free-surface flows, even in the presence of complex geometries and source terms.

\paragraph{\textbf{Related work.}} Finite element methods in space combined with 
explicit, semi-implicit, and implicit Euler discretizations in time, are a common choice for numerically solving the coupled Stokes/free-surface problem. 
The implicit Euler method allows for the largest time step sizes, 
but is computationally expensive as it requires nonlinear iterations coupling the Stokes/free-surface problem. 
On the other hand, explicit Euler and semi-implicit Euler discretizations in time \cite{JouvetPicassoRappazBlatter2008, ElmerDescrip,cheng_gong_von_neumann,WirbelJarosch,Kaus_fssa} do not require additional nonlinear iterations. 
However, these approaches lead to severe restrictions on the maximum stable time step size. 
The Free-Surface Stabilization Algorithm (FSSA) 
introduced in \cite{Kaus_fssa} is an approach for improving the maximum time step size of the semi-implicit Euler discretization in time. 
\igor{The} FSSA adds an \igorrr{nonsymmetric} surface integral term to the weak form Stokes problem to bring the full discretization closer to an implicit method. 
The approach was first introduced in the context of the mantle convection applications \cite{Kaus_fssa} and was later extended to ice sheet dynamics applications in \cite{lofgren_fssa1,lofgren_fssa2,tominec2024weakformshallowice}. 
While the authors report significant time step size improvements, the time step restriction is 
not a priori known and the overall theoretical understanding of FSSA is still missing. 
\igorrr{In \cite{Martinez}, the FSSA stabilization terms were modified to preserve the symmetry of the linear systems.}
\igorrr{In \cite{rose_buffet_heister_u_dot_n}, the authors 
give a procedure for computing the largest stable time-step for the Stokes/free-surface system 
simplified to an associated linear homogeneous eigenvalue problem, 
by modifying the FSSA stabilization terms \cite{Kaus_fssa}.}
\begin{figure}[h!]
  \centering
\begin{tabular}{cc}
  \textbf{Domain} & \hspace{1cm}\textbf{Free-surface height function} \vspace{0.2cm} \\
  \raisebox{0cm}{\includegraphics[width=0.4\linewidth]{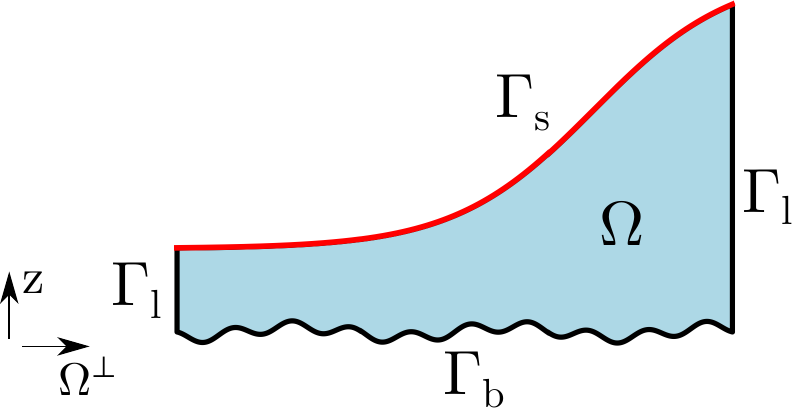}} &
  \raisebox{0.05cm}{\hspace{1cm}\includegraphics[width=0.29\linewidth]{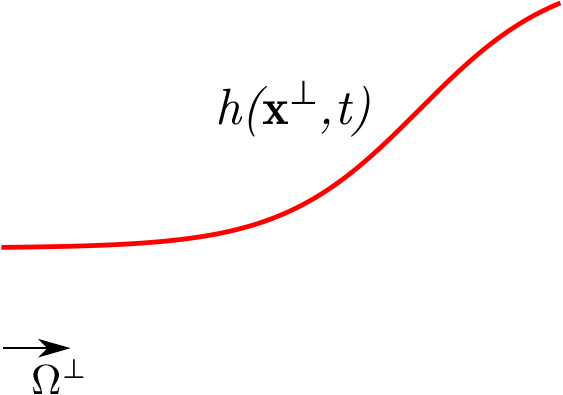}} \\
  (a) & \hspace{0.7cm} (b)
  \end{tabular}
  \caption{Plot (a) shows an example of domain $\Omega \subset \mathbb{R}^d$ with boundary split into $\Gamma_s$ (free-surface), $\Gamma_b$ (bottom topography/bedrock), 
  and $\Gamma_l$ (the lateral boundaries). 
  Plot (b) shows an explicit form free-surface height function $\surf=\surf(\bm x^\perp, t)$, where $\bm x^\perp$ is evaluated over free-surface domain $\Omega^\perp \subset \mathbb{R}^{d-1}$.}
  \label{fig:icesheet_sketch}
\end{figure}
In \cite{Audusse_FEM_Navier_stokes_Free_Surface}, the authors were the first to examine the numerical stability properties of the coupled Navier-Stokes/free-surface problem, 
discretized using the finite element method and the explicit Euler method. 
To improve the maximum time step size, the authors supplemented the free-surface equation and the Navier-Stokes problem with additional first-order consistent stabilization terms. 
Their approach improves the time step \igor{sizes up to a CFL-type} restriction that can be evaluated up to an unknown constant. \igor{Furthermore, 
the authors consider Newtonian fluids without a source term in the free-surface equation.}
In \cite{bueler2024_surf_elevation_errors}, the author gives surface height error estimates for the implicit Euler discretized coupled Stokes/free-surface problem for non-Newtonian fluids, 
however, the stability properties of computationally less expensive approaches were not investigated. 
Moreover, the estimate coupling the momentum energy balance and the surface height energy through a total energy norm is missing. 

\paragraph{\textbf{Contributions of the present paper.}} In this paper, we present a novel theoretical study of the stability 
of the coupled Stokes/free-surface problem for Newtonian and non-Newtonian fluids with a source term, in both the continuous and fully discrete settings. 
On the continuous level, 
we derive a total energy balance that couples the energy of the Stokes problem with that of the free-surface equation, 
\igor{which is key to understanding the stability.} \reviewerOne{We have not found such a total energy balance in the literature, and it is a novel contribution of this paper.}
From the total energy balance, we establish \reviewerTwo{an $L^2$-norm total energy} stability estimate and \reviewerTwo{additionally} the conservation of domain volume.
On the discrete level, we show that the implicit Euler method is unconditionally stable \reviewerTwo{when free-surface height is Lipschitz continuous at all times} and preserves domain volume. 
Furthermore, we provide numerical counterexamples showing that the semi-implicit Euler method does not conserve domain volume. 
Moreover, the FSSA from \cite{Kaus_fssa} only partially mitigates the stability issues of both the explicit and the semi-implicit Euler schemes, 
resulting in a conditional time step constraint. 
\reviewerOne{
Issues with instabilities due to a too large time step are a known problem by numerical experiments, documented for instance in \cite{lofgren_fssa1,lofgren_fssa2,tominec2024weakformshallowice}. 
However, here we prove this observation theoretically.} 
To address this limitation, we use our stability estimates to derive an optimally scaled stabilization term 
that penalizes the normal velocity at the free surface in the Stokes problem. 
This \reviewerOne{new} stabilization ensures that the explicit Euler method used for advancing the free-surface equation is unconditionally stable with respect to the time step size 
\reviewerTwo{when free-surface height is Lipschitz continuous at all times}. 
Our formulation also extends to cases with a source term in the free-surface equation. 
\reviewerOne{While stabilizations of this type are proposed in \cite{Kaus_fssa,lofgren_fssa1,lofgren_fssa2}, none of them is theoretically analyzed 
in terms of energy estimates and none of them enables arbitrary time step sizes. The first-order viscosity stabilization of the free-surface equation proposed 
in \cite{Audusse_FEM_Navier_stokes_Free_Surface} improves stability, but does not lead to arbitrary time step sizes. 
Our new stabilization is the first one that theoretically guarantees arbitrary time step sizes for the explicit Euler discretization of the coupled Stokes/free-surface problem.}
In addition, we prove that our new stabilized explicit scheme preserves the volume of the domain and numerically verify that it maintains first-order accuracy under simultaneous refinement of the mesh and the time step. 
Our analysis does not rely on any simplifications of the geometry or the original Stokes/free-surface formulation and is valid in both two and three spatial dimensions.

\paragraph{\textbf{Paper organization.}} This paper is structured as follows. 
In Section \ref{sec:model_formulation} 
we state the Stokes problem coupled to the free-surface equation in strong and weak 
form and define the domains and their boundary representations.
In Section \ref{sec:continuous_problem_analysis} and Section \ref{sec:discrete_problem_analysis} 
we present the continuous and discrete stability analyses, respectively. 
In Section \ref{sec:experiments} we provide numerical experiments verifying our theory for the Newtonian fluid case and non-Newtonian fluid case.
In Section \ref{sec:final_remarks} we give final remarks.

\section{Mathematical model}
\label{sec:model_formulation}
\subsection{Governing equations}
\label{sec:model_strong_form}
The fluid is represented as a \reviewerOne{Lipschitz domain} $\Omega = \Omega(t) \subset \mathbb{R}^d$ that evolves in time $t \in T := \lft[ 0, \hat{t} \, \rgt]$. 
The evolving free boundary (one part of $\partial\Omega$) is represented as the surface height $\surf ( \bm x^\perp,t )$ (\reviewerOne{function graph}) with 
$\bm x^\perp$ evaluated over a \igor{horizontal (flat)} \reviewerOne{Lipschitz} domain $\Omega^\perp \subset \mathbb{R}^{d-1}$. 
The remaining part of $\partial \Omega$ is a \reviewerTwo{static topography/bedrock $b(\bx^\perp ) \in C^{0,1}(\Omega^\perp)$} (\reviewerOne{function graph}). 
The area between $\surf$ and $b$ is the fluid thickness $\surf - b > 0$.
The relation between $\Omega$, $\surf$, \reviewerTwo{$b$} 
and $\Omega^\perp$ is illustrated in Figure \ref{fig:icesheet_sketch}, 
\reviewerTwo{whereas detailed mathematical definitions are given later in Section \ref{sec:domain_representation}.} 
The coupled Stokes/free-surface problem reads \reviewerTwo{as follows}. Find the fluid surface height $\surf: \Omega^\perp \times T \to \mathbb{R}$, 
velocity $\bm u: \Omega(t)  \to \mathbb{R}^d$, and pressure $\pi: \Omega(t) \to \mathbb{R}$, such that: 
\begin{equation}
  \label{eq:stokes_free_surface_strongform}
  \left\{
    \begin{aligned}
      \partial_t \surf + \bm u^\perp \cdot \nabla^\perp \surf - u_z &= a,\\
      - \nabla \cdot  ( 2\, \mu(\bm u) \bm D \bm u) + \nabla \pi &= - \rho g \hat{\bz},\\
\nabla \cdot  \bm u &= 0,\\
    \end{aligned}
  \right.
\end{equation}
where
$\bm u^\perp: \Omega^\perp \times T \to \mathbb{R}^{d-1}$ is a vector of the horizontal surface velocities, 
$\nabla^\perp$ is the \igor{horizontal} gradient operator over $\Omega^\perp$, 
$u_z$ is the vertical surface velocity, 
$a: \Omega^\perp \times T \to \mathbb{R}$ is the source term with regularity \reviewerTwo{$a \in L^2([0,\hat t];\, L^2(\Omega^\perp))$}, $\nabla$ is the gradient operator over $\Omega$, $\bm D \bm u  = \frac{1}{2} \left(  \nabla \bm u + (\nabla \bm u)^T \right) $ is the strain rate tensor, $g>0$ is the gravitational acceleration, $\rho>0$ is the fluid density, and $\mu: \Omega \to \mathbb{R}^+$ is a viscosity function. 
We consider a family of viscosity functions in the power-law form:
\begin{equation}
  \label{eq:stokes_viscosity}
\mu(\bm u) = \mu_0\,|\bm D \bm u|^{p-2},\quad 1< p \leq 2,
\end{equation}
where $\mu_0 \in \mathbb{R}^+$ is an effective viscosity constant, $|\bm D \bm u| = (\sum_{i,j=1}^d (\bm D \bm u)_{ij}\, (\bm D \bm u)_{ij})^{1/2}$ is the Frobenius norm of $\bD \bu$. 
\igor{Glaciologists} use Glen's flow law exponent $n$ in place of $p$, where the relation is $p=(1/n)+1$ \cite{Glen1955,Hirnpowerlaw}. 
When $p=2$, the fluid is Newtonian. When $1<p<2$, the fluid is non-Newtonian.
We supplement \eqref{eq:stokes_free_surface_strongform} with an initial condition for the surface $\surf \lft(\bx^\perp, 0 \rgt) = \surf_0 \lft( \bx^\perp \rgt)$, 
and the boundary conditions:
\begin{equation}
  \label{eq:stokes_bcs}
    \begin{aligned}
      \bu = \bm 0 \text{ on } \Gamma_b,\quad
      \bm u \cdot \bm n= 0 \text{ on } \Gamma_{l}, \quad
      \bm \sigma \cdot \bm n  = \bm 0 \text{ on } \Gamma_s, \quad 
      \reviewerTwo{(\bm \sigma \cdot \bm n)_t = \bm 0 \text{ on } \Gamma_l,}
    \end{aligned}
\end{equation}
where $\Gamma_b$ is the bedrock, $\Gamma_s$ is the free-surface, $\Gamma_l$ are the lateral boundaries, $\bm \sigma=2 \mu(\bu) \bD  \bu  - \pi \polI $ 
is the Cauchy stress tensor, \reviewerTwo{and $(\bm \sigma \cdot \bm n)_t$ denotes the tangential component of the normal stress tensor.}

\reviewerOne{
  The mathematical model as given in \eqref{eq:stokes_free_surface_strongform}, 
  \eqref{eq:stokes_viscosity}, and \eqref{eq:stokes_bcs},  provides a satisfactory description 
  of free-surface flows when the advective forces are negligible compared to viscous forces. } \reviewerTwo{In the thin-domain regime relevant to ice sheets, diffusion dominates under no-slip basal boundary conditions, and the fully coupled system is well approximated by a parabolic equation \citep{greve_book}.}
  \reviewerOne{ Moreover, the model is applicable to incompressible Newtonian fluids, as well as incompressible 
  non-Newtonian fluids with shear-thinning behavior. 
  Applications are found for instance in ice sheet dynamics \cite{greve_book}, mantle convection \cite{Kaus_fssa}, 
  and lava flow \cite{2020_hinton_lava}.
}

\subsection{Weak formulation and function spaces}
\label{sec:weak_form}
We analyze the weak formulation of the coupled problem \eqref{eq:stokes_free_surface_strongform}, which is defined as follows. 
Find $\bm u \in \bm V(\Omega)$, $\pi \in Q(\Omega)$, and $\surf \in Z(\Omega^\perp) \times \calC^1(T) $ such that: 
\begin{equation}
  \label{eq:stokes_free_surface_weakform}
  \left\{
\begin{aligned}
  (\partial_t \surf, w)_{\Omega^\perp} + \lft( \bu^\perp \cdot \nabla ^\perp \surf - u_z, w \rgt)_{\Omega^\perp}   &= (a,w)_{\Omega^\perp}, \\
  2 \lft( \sqrt{\mu}\,\bm D \bm u, \sqrt{\mu}\,\bm D \bm v  \rgt)_{\Omega(t)} - (\pi, \nabla \cdot  \bv)_{\Omega(t)} &= - \rho g(  \hat{\bz}, \bv)_{\Omega(t)},\\
( \nabla \cdot  \bu, q)_{\Omega(t)} &= 0,
\end{aligned}
\right.
\end{equation}
for all $\bm v \in \bm V(\Omega)$, $q \in Q(\Omega)$, and $w \in Z(\Omega^\perp)$, 
\reviewerTwo{where inner products and norms over time-dependent domains are for 
scalar fields (and analogously vector fields) 
$f: \Omega \times T \to \mathbb{R}$ and $g: \Omega^\perp \times T \to \mathbb{R}$ defined as:
\begin{equation}
  \label{eq:weakform:innerproduct_notation_continuous}
  \begin{aligned} 
  (f,f)_{\Omega(t)} &:=  \|f\|^2_{L^2(\Omega(t))} = \int_{\Omega(t)} f^2\, d\bx = \int_{\Omega(t)} (f(\bm x, t))^2\ d\bx& \\
  (g,g)_{\Omega^\perp} &:= \|g\|^2_{L^2(\Omega^\perp)} = \int_{\Omega^\perp} g^2\, d\bx^\perp = \int_{\Omega^\perp} (g(\bm x^\perp, t))^2\, d\bx^\perp&\quad t \in T.
  \end{aligned}
\end{equation}
Surface integrals over $\Gamma_s(t)$ are interpreted analogously as:
\begin{equation}
\int_{\Gamma_s(t)} f\, ds = \int_{\Gamma_s(t)} f(\bm x, t)\, ds,\quad t \in T.
\end{equation}
}
We now define the solution spaces $\bm V$, $Q$, and $Z$.
The Lebesgue space $L^\alpha(\Omega)$ is defined as \cite{guermond_voli}:
$L^\alpha(\Omega) = \{ q: \Omega \to \mathbb{R}\, |\, q \text{ measurable and } \|q\|^\alpha_{L^\alpha(\Omega)} < \infty \},\, \alpha \in [1, \infty),$
where an associated $L^\alpha$-norm is $\|q\|_{L^\alpha(\Omega)} = (\int_\Omega |q|^\alpha\, d \bx)^{1/\alpha}$. The
Sobolev space $\bm W^{1,\alpha}(\Omega)$ and its restriction to vanishing functions over $\Gamma_b \subset \partial\Omega$ are respectively 
defined through:
\begin{equation*}
    \begin{aligned}
\bm W^{1,\alpha}(\Omega) &= \Big\{ \bm v \in [L^\alpha(\Omega)]^d \, \Big|\, \|\bm v\|_{L^\alpha(\Omega)}^\alpha + \|\nabla \bm v\|_{L^\alpha(\Omega)}^\alpha < \infty \Big\}, \\
\bm W^{1,\alpha}_{\Gamma_b,\reviewerTwo{\Gamma_l}}(\Omega) &= \Big\{ \bm v \in \bm W^{1,\alpha}(\Omega) \, \Big|\, \bm v|_{\Gamma_b} = 0,\, \reviewerTwo{(\bm v \cdot \bm n)|_{\Gamma_l} = 0} \Big\}.
    \end{aligned}
\end{equation*}
The solution to the nonlinear Stokes problem in weak form belongs to the velocity space $\bm V(\Omega)$ and the pressure space $Q(\Omega)$, defined as \cite{Belenki2012}:
\begin{equation}
\bm V(\Omega) = \bm W^{1,p}_{\Gamma_b,\reviewerTwo{\Gamma_l}}(\Omega),\quad Q(\Omega) = L^{p'}(\Omega),\quad p' = \frac{p}{p-1},\quad 1 < p \leq 2,
\end{equation}
where $p$ is a parameter in \eqref{eq:stokes_viscosity}. 
\reviewerTwo{Function space for the free-surface height is chosen as:
$$Z ( \Omega^\perp ) \subset L^2 ( \Omega^\perp ),$$
where an $L^2$ inner product is $(q,w)_{\Omega^\perp} = \int_{\Omega^\perp} q\, w\, d \bx^\perp$. For well-posedness, a tighter function space 
is required for the surface height, see discussion in Section \ref{sec:continuous_problem_analysis:well_posedness_discussion}. 
However, explicitly defining the tighter function space $Z(\Omega^\perp)$ requires modifications of the 
standard free-surface equation which is an active research question \cite{bueler2024_surf_elevation_errors}, and we leave it for future work.
}


We approximate the solutions to the Stokes problem and the free-surface equation using the finite element method. 
The discrete domains $\Omega_h$ and $\Omega_h^\perp$ are represented by \reviewerOne{shape-regular and quasi-uniform} meshes $\mathcal T_h$ and $\mathcal T_h^\perp$. 
Each mesh is composed of a finite number of disjoint elements $K$ and $K^\perp$ so that $\Omega_h = \cup_{K\in\mathcal{T}_h} K$ and $\Omega_h^\perp = \cup_{K\in\mathcal{T}_h^\perp} K^\perp$. 
We denote $\mathcal{E}_h$ and $\mathcal{E}_h^\perp$ a set of  interior $d-1$ dimensional faces in $\mathcal{T}_h$ and $\mathcal{T}_h^\perp$ respectively. 
The associated finite element spaces are the lowest-order Taylor-Hood velocity and pressure spaces, and a piecewise linear space for free-surface height:
\begin{equation*}
    \begin{aligned}
\bm V_h(\Omega_h) &= \Big\{  \bv_h \in \mathcal{C}^0(\Omega_h) \, \Big|\, \bm v_h \big |_{K} \in \mathbb{P}_{2}(K),\, \forall K \in \mathcal T_h \Big\}, \\
Q_h(\Omega_h) &= \Big\{ q_h \in \mathcal{C}^0(\Omega_h) \, \Big|\, q_h \big |_{K} \in \mathbb{P}_1(K),\, \forall K \in \mathcal T_h \Big\}, \\
Z_h(\Omega_h^\perp) &= \Big\{ w_h \in \mathcal{C}^0(\Omega_h^\perp) \, \Big|\, w_h \big |_{K^\perp} \in \mathbb{P}_1(K^\perp),\, \forall K^\perp \in \mathcal T_h^\perp \Big\}, \\
    \end{aligned}
\end{equation*}
where $\mathbb{P}_1$, $\mathbb{P}_2$ are multivariate polynomial spaces of order $1$ and $2$. The analysis naturally extends to higher-order Taylor-Hood elements and higher-order mesh representations. 

\subsection{Surface representations and integration domains}
\label{sec:domain_representation}
In this section, we describe how the domain, surface integrals, and normal vectors depend on the surface height $\surf$, as this will be required in the subsequent analysis. For $d=2$ we define $\bx = (x,z), \bx^\perp = x$, whereas for $d=3$ we define $\bx = (x,y,z), \bx^\perp = (x,y)$. 
Then, the domain  and its boundary are defined as:
\begin{equation}
  \label{eq:surface_representations:Omega}
\Omega(t) = \lft \{\bx \in \mathbb{R}^d\,|\, \bx^\perp \in \Omega^\perp\, \text{ and }\, b(\bx^\perp) < z < \surf(\bx^\perp, t) \rgt \}, \quad \partial\Omega(t) = \Gamma_s(t) \cup \Gamma_b \cup \Gamma_l,
\end{equation}
where:
\begin{equation}
  \label{eq:surface_representations:Gamma_s_and_Gamma_b}
  \begin{aligned}
\Gamma_s(t) &= \lft \{ \lft( \bx^\perp,\surf(\bx^\perp, t) \rgt)  \in \mathbb{R}^d\,\big |\, \bx^\perp \in \Omega^\perp \rgt \},\quad \Gamma_b = \lft \{ \lft (  \bx^\perp, b(\bx^\perp) \rgt) \in \mathbb{R}^d\, \big |\, \bx^\perp \in \Omega^\perp \rgt \}, \\
\Gamma_l &= \p \Omega(t) \backslash (\Gamma_s(t) \cup \Gamma_b).
  \end{aligned}
\end{equation}
The volume of the domain is then defined as $|\Omega(t)| = \int_{\Omega(t)} d \bx = \int_{\Omega^\perp} (\surf(\bm x,t) - b(\bm x)) d \bx^\perp$.

When computing integrals over $\Gamma_s(t)$, it is convenient to change the integration domain to the time-independent evaluation domain $\Omega^\perp$. 
Let $\bm r_s = (\bx^\perp,\surf(\bx^\perp,t)) = (x, y, \surf(x,y,t))$ 
be a parametrization of $\Gamma_s(t)$. Then we have:
\begin{equation}
  \label{eq:surface:integration_domains}
  \int_{\Gamma_s(t)} f(\bm x)\, ds = \int_{\Omega^\perp} f(\bm r_s)\, \|\partial_{x} \bm r_s \times \partial_y \bm r_s \|_2\, d \bx^\perp 
  = \int_{\Omega^\perp} f(\bm r_s)\, \|(-\nabla^\perp \surf, 1)\|_2 \, d \bx^\perp,
\end{equation} 
where $\|\cdot\|_2$ denotes the Euclidean norm. This follows from the fact that 
$\partial_x \bm r_s = (1, 0, \partial_x \surf)$, $\partial_y \bm r_s = (0, 1, \partial_y \surf)$, 
and thus $\partial_x \bm r_s \times \partial_y \bm r_s = (-\partial_x \surf, -\partial_y \surf, 1) = (-\nabla^\perp \surf, 1)$. 
Using that relation, the outward unit normal vector $\bm n$ to $\Gamma_s$ is:
\begin{equation}
  \label{eq:surface:normals}
  \bm n(\bm x, t) = \frac{\partial_x \bm r_s \times \partial_y \bm r_s}{\|\partial_{x} \bm r_s \times \partial_y \bm r_s \|_2}  
  = \frac{(-\nabla^\perp \surf, 1)}{\|(-\nabla^\perp \surf, 1)\|_2} , \quad \bm x \in \Gamma_s(t) .
\end{equation}
Finally, the surface velocity $\bu^\perp : \Omega^\perp \times T \to \mR^{d-1} $ is denoted as
  $\bu^\perp =  \lft( u_x \lft(\bx^\perp, \surf\rgt), u_y \lft(\bx^\perp, \surf \rgt) \rgt)$ when $d=3$ and $\bu^\perp =  \lft( u_x \lft( \bx^\perp,\surf \rgt) \rgt)$ when $d=2$.
Note that the relations $  \int_{\Gamma_s} f(\bm x)\, ds = \int_{\Omega^\perp} f(\bm r_s)\, \|(-\nabla^\perp \surf, 1)\|_2 \, d \bx^\perp$ 
and $\bm n = \frac{(-\nabla^\perp \surf, 1)}{\|(-\nabla^\perp \surf, 1)\|_2}$ derived above also applies for $d=2$.

\section{Continuous problem analysis}
\label{sec:continuous_problem_analysis}
In this section, we analyze the stability of the problem as stated in \eqref{eq:stokes_free_surface_weakform}, posed over the infinite-dimensional spaces $\bm V(\Omega)$, $Q(\Omega)$, and $Z(\Omega^\perp)$.
Our goal is to derive the total energy balance, which couples the energy balance of the momentum equations 
to the free-surface height energy. 
The framework developed here is also used in later sections to analyze stability and conservation in the discrete setting, where explicit, semi-implicit, and implicit Euler methods are used for time discretization, and finite elements are used in space. 

\subsection{\reviewerOne{Discussion on well-posedness using $L^2$-norm estimates as a building block}}
\label{sec:continuous_problem_analysis:well_posedness_discussion}
\reviewerOne{In this section, we discuss the requirements for well-posedness of the coupled Stokes/free-surface problem in \eqref{eq:stokes_free_surface_weakform} and the importance of $L^2$-stability of the free surface $\surf$.

\paragraph{Nonlinear Stokes well-posedness requirements.} Well-posedness of the nonlinear Stokes problem in \eqref{eq:stokes_free_surface_weakform} 
relies on coercivity of the stress term $(\sqrt{\mu}\,\bm D \bm u, \sqrt{\mu}\,\bm D \bm v)$, 
the inf-sup condition for $(\pi, \nabla\cdot \bm v)$, and continuity of both forms \cite{Belenki2012}.
The inf-sup condition on $V(\Omega)\times Q(\Omega)$ follows from the Nečas inequality \cite{Girault_and_Raviart}.
To address coercivity on $V(\Omega)$, we recall the relation between $(\sqrt{\mu}\,\bm D \bm u,\sqrt{\mu}\,\bm D \bm u)_{L^2(\Omega)}$
and $\|\bm D \bm u\|_{L^p(\Omega)}$. Using viscosity as defined in \eqref{eq:stokes_viscosity}, we have:
\begin{equation}
\label{eq:L_p_to_L_2_strain}
(\sqrt{\mu}\,\bm D \bm u, \sqrt{\mu}\,\bm D \bm u)_{L^2(\Omega)}
= \int_\Omega \mu\,|\bm D \bm u|^2\,d\bx 
= \int_\Omega \mu_0\,|\bm D \bm u|^{p-2}\, |\bm D \bm u|^{2}\,d\bx
= \int_\Omega \mu_0\,|\bm D \bm u|^{p}\,d\bx
= \mu_0\,\|\bm D \bm u\|_{L^p(\Omega)}^{p}.
\end{equation}
Applying Korn's inequality and the Poincaré inequality gives the coercivity estimate
\cite{Belenki2012}: $(\sqrt{\mu}\,\bm D \bm u, \sqrt{\mu}\,\bm D \bm u)_{L^2(\Omega)}
\geq C\,\mu_0\,\|\bm u\|_{W^{1,p}(\Omega)}^{p}.$ A sufficient assumption for both coercivity and the inf-sup condition
is Lipschitz continuity of the domain  $\partial\Omega\in C^{0,1}$ \cite{Belenki2012,Girault_and_Raviart}.

\paragraph{Free-surface equation regularity.} The free-surface equation in \eqref{eq:stokes_free_surface_weakform} is on its own (i.e. viewed as de-coupled from the Stokes equation)  
a linear advection equation with source terms.
In the present setting, the free surface is described by a height function, function graph $\surf$, defined over $\Omega^\perp$.
An $L^2(\Omega^\perp)$ estimate of $\surf$ does not guarantee that the graph of $\surf$ defines a sufficiently regular boundary for the Stokes problem. 
Moreover, as $\partial\Omega=\partial\Omega(\surf)$, the Stokes requirement of the 
Lipschitz continuity of the domain boundary requires $\surf$
to also be Lipschitz continuous on $\Omega^\perp$, which is a requirement going beyond the basic $L^2$-control. 
However, the regularity of $\surf$ may be increased by adding regularization terms to the free-surface equation \cite{bueler2024_surf_elevation_errors}.
If such a term induces a semi-norm, $L^2$-estimates of $\surf$ may be combined with it to obtain a stronger norm control.
This motivates the role of $L^2$-stability of the free surface, which we view as a fundamental building block 
for well-posedness of the coupled problem.

\paragraph{Discussion on improving regularity of $\surf$.} 
We now examine possible regularity improvements using Sobolev embeddings \cite[Theorem 4.12]{Adams_Fournier_book}. 
The space $L^2(\Omega^\perp)$ does not embed into $C^{0,1}(\Omega^\perp)$, 
and thus cannot guarantee Lipschitz continuity of $\partial\Omega$. In physically relevant settings, $\Omega$ is three-dimensional,
so $\surf$ is defined on a two-dimensional $\Omega^\perp$. In this case, even $W^{1,2}(\Omega^\perp)$ does not embed into $C^{0,1}(\Omega^\perp)$,
and adding a standard first order viscosity term controlling $W^{1,2}$-norm to the free-surface equation is insufficient.
In \cite{bueler2024_surf_elevation_errors}, a $W^{1,4}(\Omega^\perp)$ control term is considered, 
but the embedding $W^{1,4}(\Omega^\perp)\hookrightarrow C^{0,1/2}(\Omega^\perp)$ still falls short of Lipschitz continuity when $\Omega^\perp \subset \mathbb{R}^2$. 
Furthermore, the term introduces additional nonlinearity into the free-surface equation which increases the computational 
complexity when solved numerically.
An alternative is to exploit the embedding $W^{2,2}(\Omega^\perp)\hookrightarrow C^{0,\alpha_L}(\Omega^\perp)$ with $\alpha_L<1$, 
which also does not guarantee Lipschitz continuity, but it gets infinitesimally close to it. This regularity is possible to achieve 
through a linear 
$W^{2,2}$ (hyperviscosity) stabilization term, for instance via interior penalty approach for Lagrange-type finite elements \cite{Brenner_2005}. 
Full Lipschitz continuity follows for instance from $W^{2,2+\epsilon}(\Omega^\perp)\hookrightarrow C^{0,1}(\Omega^\perp)$ for any $\epsilon>0$, but controlling the corresponding norm of $\surf$ 
likely requires adding nonlinear terms to the free-surface equation. Another 
approach is exploiting the embedding $W^{3,2}(\Omega^\perp)\hookrightarrow C^{1}(\Omega^\perp)$, but controlling the corresponding norm of $\surf$
is considerably more challenging in finite element discretizations that are typically only conforming with $W^{1,2}$ spaces. 
When choosing an embedding and the corresponding regularization term, it is also important to make the regularization physically meaningful and consistent with the original model.

\paragraph{Regularization summary.} Regularizing the free-surface equation by control terms that ensure sufficient regularity of $\surf$,
ideally Lipschitz continuity, is a prerequisite for well-posedness of the coupled problem \eqref{eq:stokes_free_surface_weakform}. 
In this paper, we do not focus on the regularization of the surface continuity from a theoretical perspective. 
Instead, we adopt a pragmatic approach and introduce a high-order, hyperviscosity-type ($H^2$-like) stabilization term
in the discretized free-surface equation. This term is treated implicitly in time
to avoid additional time-step restrictions. A rigorous continuous well-posedness analysis of this regularization
is beyond the scope of the present work and will be addressed in future work. 
Instead, we give $L^2$-norm estimates which are a fundamental building block for well-posedness when the free-surface equation 
is regularized for Lipschitz continuity, as motivated throughout this 
section. The new framework developed in this paper can directly be used when analyzing well-posedness after the free-surface equation 
is appropriately regularized.
}



\subsection{Relation between the Stokes problem on $\Omega$ and the free-surface equation on $\Omega^\perp$}
We establish two relations used for the stability analysis of the coupled Stokes/free-surface problem in the sections that follow. 
The first relates one term in the free-surface equation on $\Omega^\perp$ to an equivalent formulation on $\Gamma_s$. 
The second relates the term formulated on $\Gamma_s$, to the right-hand side of the Stokes equation. 
\begin{lemma}[Free-surface equation to surface integral]
\label{lemma:free_surf_equation_to_surface_integral}
Let $\surf: \Omega^\perp \times T \to \mathbb{R}$ be the surface height function associated with 
$\Gamma_s(t)$ \eqref{eq:surface_representations:Gamma_s_and_Gamma_b}. 
\reviewerTwo{Assume $\surf(t) > b$.} Then the following relation holds: 
\begin{equation}
  \label{eq:free_surf_eq_to_surface_integral}
  (- \bu^\perp \cdot \nabla ^\perp \surf + u_z, w)_{\Omega^\perp} = \int_{\Gamma_s(t)} (\bm u \cdot \bm n)\, w\, ds,\quad \forall w \in Z (\Omega^\perp),
\end{equation}
and when $w = \surf$, the relation becomes
  $(- \bu^\perp \cdot \nabla ^\perp \surf + u_z, \surf )_{\Omega^\perp} = \int_{\Gamma_s(t)} (\bm u \cdot \bm n)\, z\, ds$ with 
  $(x,y,z) \in \Omega(t)$.
\end{lemma}
\begin{proof}
We rewrite the integration domain from $\Omega^\perp$ to $\Gamma_s(t)$. We have 
  $      (- \bu^\perp \cdot \nabla ^\perp \surf + u_z, w)_{\Omega^\perp} = \int_{\Omega^\perp} (- \bu^\perp \cdot \nabla ^\perp \surf + u_z)\,w\, d \bx^\perp   = \int_{\Omega^\perp}  \bu \cdot (-\nabla^\perp \surf, 1)\, w\, d \bx^\perp =  \int_{\Omega^\perp}  \bu \cdot ( -\nabla^\perp \surf, 1)\, w\,   \frac{\| ( -\nabla^\perp \surf, 1) \|_2 }{\| ( -\nabla^\perp \surf, 1) \|_2} d \bx^\perp$.
  It follows that: (i) the vector $\frac{ ( -\nabla^\perp \surf, 1)  }{\| ( -\nabla^\perp \surf, 1) \|_2}$ corresponds to the unit normal $\bm n$ as defined in \eqref{eq:surface:normals}, and (ii) the remaining factor
  $\| ( -\nabla^\perp \surf, 1) \|_2$ corresponds to the arclength $ds$ of the surface $\Gamma_s(t)$, as defined in \eqref{eq:surface:integration_domains}. Using that, we rewrite the original integral as:
  \begin{equation} \label{eq:lemma_free_surface}
    (- \bu^\perp \cdot \nabla ^\perp \surf + u_z, w)_{\Omega^\perp} = \int_{\Omega^\perp}  \bm u \cdot \frac{( -\nabla^\perp \surf, 1)}{\| ( -\nabla^\perp \surf, 1) \|_2}\, w\,\| ( -\nabla^\perp \surf, 1) \|_2\, d \bx^\perp = \int_{\Gamma_s(t)} (\bm u \cdot \bm n)\, w\, ds,
  \end{equation}
which proves the first statement \eqref{eq:free_surf_eq_to_surface_integral} of the lemma.
In \eqref{eq:lemma_free_surface} we set $w=\surf$. 
Since $\Gamma_s(t)$ is in \eqref{eq:surface_representations:Gamma_s_and_Gamma_b} parameterized as $z = \surf(\bx^\perp, t)$ and since we are integrating over $\Gamma_s$, we substitute $\surf=z$, 
which proves the second statement of this proposition.
\end{proof}
\begin{lemma}[Stokes problem to surface integral]
  \label{lemma:Stokes_to_surface_integral_ztrick}
  Let $\bm u \in \bm V(\Omega)$ solve the Stokes problem in \eqref{eq:stokes_free_surface_weakform}. 
Then, the relation between the volume integral $- \rho g ( \hat{\bz}, \bm u)_{\Omega(t)}$ and the surface integral over $\Gamma_s(t)$ is: 
\begin{equation}
  - \rho g( \hat{\bz}, \bm u)_{\Omega(t)} = - \rho g \, \int_{\Gamma_s(t)} (\bm u \cdot \bm n)\, z\, ds,\quad \bx  \in \Omega(t).
\end{equation}
\end{lemma}
\begin{proof}
  We write the gravity term as $- \rho g \hat{\bz} = - \rho g \nabla z$. It follows that:
  $$- \rho g ( \hat{\bz}, \bm u)_{\Omega(t)} = - \rho g (\nabla z, \bm u)_{\Omega(t)} 
  = \rho g (\nabla \cdot \bm u, z)_{\Omega(t)} - \rho g \int_{\partial\Omega(t)} (\bm u \cdot \bm n)\, z\, ds = - \rho g \int_{\Gamma_s(t)} (\bm u \cdot \bm n)\, z\, ds.$$
In the first step, we integrated $(\nabla z, \bm u)_{\Omega(t)}$ by parts. In the second step, we used the incompressibility condition $(\nabla \cdot \bm u, z)_{\Omega(t)} = 0$ from \eqref{eq:stokes_free_surface_weakform}, noting that $z \in L^2(\Omega)$. 
Finally, the boundary integral was restricted to $\Gamma_s$ using the boundary conditions \reviewerTwo{on $V(\Omega)$}.
\end{proof}

\subsection{Stability analysis}
In this section, we show the stability of the coupled Stokes/free-surface problem \eqref{eq:stokes_free_surface_weakform} in the continuous setting, 
through an energy balance and a corresponding total energy estimate, \igor{which is the key to understanding the stability also in discretized settings}. In the later sections, we use this total energy estimate 
to argue about the stability of the numerical schemes.
\begin{theorem}[Continuous \reviewerTwo{$L^2$} stability]
  \label{theorem:surface_height_energy_rate}
  Let the surface height $\surf$ \reviewerTwo{(with $\surf>b$)} and the velocity $\bm u$ solve the coupled system \eqref{eq:stokes_free_surface_weakform} on domain $\Omega = \Omega(t)$, with source term $a=a \lft(\bx^\perp, t \rgt)$,
  viscosity $\mu=\mu(\bm u)$, fluid density $\rho$ and the gravitational acceleration $g$.
  Then, the following total energy balance holds:
  \begin{equation} \label{eq:theorem_energy_balance}
  \frac{1}{2} \partial_t \|\surf\|_{L^2(\Omega^\perp)}^2 = - \frac{2}{ \rho g}\, \|\sqrt{\mu}\,\bm D \bm u\|^2_{L^2(\Omega(t))} + (a,\surf)_{\Omega^\perp}.
  \end{equation}
  Moreover, the following strong stability estimate in the total energy norm holds:
  \begin{equation} \label{eq:theorem_stability}
    \begin{split}
      \|\surf(t)\|_{L^2(\Omega^\perp)}^2 + \frac{4}{\rho g}\, \int_0^{\, t} \|\sqrt{\mu}\,\bm D \bm u\|^2_{L^2(\Omega(\tau))} d\tau &\leq \|\surf_0\|_{L^2(\Omega^\perp)}^2 
      + 2 \reviewerTwo{\|\surf_0 \|_{L^2(\Omega^\perp)}}\, \int_0^{\, t} \| a \|_{L^2(\Omega^\perp)}\,  d\tau \\
       & \quad 
       + \reviewerTwo{2 \left(\int_0^{\, t} \| a \|_{L^2(\Omega^\perp)}  d\tau\right)^2} \qquad \forall t \in T.
    \end{split}
  \end{equation}
  \end{theorem}
  \begin{proof}
    \,\\
\textbf{(i) Energy balance.}  
Setting $\bm v = \bm u$ and $q = \pi$ in \eqref{eq:stokes_free_surface_weakform}, yields 
$2\,\|\sqrt{\mu}\, \bm D \bm u\|^2_{L^2(\Omega(t))}= - \rho g ( \hat{\bz}, \bm u)_{L^2(\Omega(t))}$. 
Applying Lemma \ref{lemma:Stokes_to_surface_integral_ztrick} to the right-hand side gives $2\,\|\sqrt{\mu}\,\bm D \bm u\|_{L^2(\Omega(t))}^2 = - \rho g \int_{\Gamma_s(t)} (\bm u \cdot \bm n)\, z\, ds.$ We now: (i) apply 
     Lemma \ref{lemma:free_surf_equation_to_surface_integral} to rewrite the surface integral over $\Gamma_s(t)$ as an integral 
     over $\Omega^\perp$, and (ii) use free-surface equation in \eqref{eq:stokes_free_surface_weakform} with $w=\surf$ to obtain
     \begin{equation}
      \label{eq:lemma:surface_height_energy_rate_tmp1}
      \begin{aligned}
     2\,\|\sqrt{\mu}\,\bm D \bm u\|_{L^2(\Omega(t))}^2 &= - \rho g \, (- \bu^\perp \cdot \nabla ^\perp \surf + u_z, \surf)_{\Omega^\perp} = - \rho g (\partial_t \surf, \surf)_{\Omega^\perp} + \rho g  (a,\surf)_{\Omega^\perp}  \\
     & = -\frac{ \rho g}{2} \partial_t \|\surf\|^2_{L^2(\Omega^\perp)} + \rho g (a,\surf)_{\Omega^\perp}. 
      \end{aligned}
     \end{equation}
Next, we divide both sides of \eqref{eq:lemma:surface_height_energy_rate_tmp1} by $ \rho g$. Rearranging the terms then proves the first statement  \eqref{eq:theorem_energy_balance} of this proposition.

\noindent
\textbf{(ii) Stability estimate.}  
Using \eqref{eq:theorem_energy_balance} we have:
\begin{equation}
  \label{eq:lemma:surface_height_energy_estimate_tmp1}
\frac{1}{2} \partial_t\|\surf\|^2_{L^2(\Omega^\perp)} = \|\surf\|_{L^2(\Omega^\perp)}\,\partial_t\|\surf\|_{L^2(\Omega^\perp)} \leq - \frac{2}{ \rho g}\, \|\sqrt{\mu}\,\bm D \bm u\|^2_{L^2(\Omega(t))} + \|a\|_{L^2(\Omega^\perp)}\, \|\surf\|_{L^2(\Omega^\perp)}.
\end{equation} 
In the above relation we use $- \frac{2}{ \rho g}\, \|\sqrt{\mu}\,\bm D \bm u\|^2_{L^2(\Omega(t))} \leq 0$ and 
then divide the resulting inequality by $\|\surf\|_{L^2(\Omega^\perp)}$, which gives:
$ \partial_t \| \surf(t) \|_{L^2(\Omega^\perp)} \leq \|a\|_{L^2(\Omega^\perp)}.$ 
We integrate the latter relation in time to arrive at:
$ \| \surf(t) \|_{L^2(\Omega^\perp)} \leq \|\surf_0\|_{L^2(\Omega^\perp)} + \int_0^{\, t}\, \|a(\tau)\|_{L^2(\Omega^\perp)}\, d\tau.$  
Inserting this on the right-hand-side of \eqref{eq:lemma:surface_height_energy_estimate_tmp1}, 
and then time-integrating and multiplying by $2$ on both sides of the inequality, finalizes the proof of the second statement \eqref{eq:theorem_stability} 
of this proposition.
\end{proof}

\subsection{Volume conservation}  \label{Sec:volume_conservation}
In the continuous setting, the domain volume $|\Omega|$ is conserved over time. 
For completeness, we briefly revisit the volume conservation analysis. 
\reviewerTwo{We first assume $\surf>b$.} 
Then, using the free-surface equation \eqref{eq:stokes_free_surface_weakform} with $w=1$ we obtain:
\begin{equation}
  \label{eq:continuous_volume_conservation_tmp1}
  \begin{aligned}
\partial_t |\Omega(t)| & = \partial_t \int_{\Omega^\perp} (\surf-b)\, d \bx^\perp = \int_{\Omega^\perp} \partial_t \surf\, d \bx^\perp = (\partial_t \surf, 1)_{\Omega^\perp} = ( -\bu^\bot \cdot \nabla ^\bot \surf + u_z + a, 1)_{\Omega^\bot}.
  \end{aligned}
\end{equation}
Next, applying Lemma \ref{lemma:free_surf_equation_to_surface_integral} to relate the integration over $\Omega^\perp$ to integration over the surface $\Gamma_s(t)$, and using integration by parts together with the incompressibility condition from \eqref{eq:stokes_free_surface_weakform} with $q = 1 \in L^2(\Omega)$
we obtain:
\begin{equation} \label{eq:continuous_volume_conservation_tmp2}
\begin{aligned}
  \partial_t |\Omega(t)|  = ( -\bu^\bot \cdot \nabla ^\bot \surf + u_z + a, 1)_{\Omega^\bot}  = \int_{ \p \Omega} \bm u \cdot \bm n \, ds + \int_{\Omega^\perp} a\, d \bx^\perp &= \int_{ \Omega(t)} \!\!\!\nabla \cdot  \bm u\, d \bx + \int_{\Omega^\perp} a\, d \bx^\perp \\
  &=  \int_{\Omega^\perp} a\, d \bx^\perp\!\!.
    \end{aligned}
  \end{equation}  
This shows that the domain volume only changes due to the source term $a$.

\section{Discretized Stokes/free-surface problems and their analysis}
\label{sec:discrete_problem_analysis}
In this section, we analyze stability and volume conservation 
when the finite element method is used in space and explicit Euler, semi-implicit Euler, and implicit Euler are used to discretize time. 
We aim to derive discrete analogues of the continuous total energy estimate \eqref{eq:theorem_stability} and volume conservation balance \eqref{eq:continuous_volume_conservation_tmp2}. 
\reviewerOne{
  We state the discrete estimates as one‑step energy inequalities. When these are iterated, 
  they give a Grönwall‑type bound in terms of accumulated forcing, rather than a uniformly decaying (dissipative) energy estimate.
}
\reviewerTwo{In all of the analysis, we assume that positivity $\surf-b>0$ everywhere on $\Omega^\perp$ holds naturally without enforcing it explicitly in the numerical schemes. 
Enforcing positivity is a separate research topic, and we refer to \cite{bueler2024_surf_elevation_errors} for a discussion on this topic in the context of the coupled Stokes/free-surface problem.}
\reviewerTwo{
Throughout this section we use the notation for integration, inner products, and norms, analogously to 
\eqref{eq:weakform:innerproduct_notation_continuous}, with the time variable fixed at $t=t^n \in T$, $n=0,1,2,...$, hence: 
\begin{equation}
  \label{eq:weakform:innerproduct_notation_discrete}
  \begin{aligned} 
  (f,f)_{\Omega_h^n} &:=  \|f\|^2_{L^2(\Omega_h^n)} = \int_{\Omega_h^n} f^2\, d\bx = \int_{\Omega_h^n} (f(\bm x, t^n))^2\ d\bx,&\\
  (g^n,g^n)_{\Omega_h^\perp} &:= \|g^n\|^2_{L^2(\Omega_h^\perp)} = \int_{\Omega_h^\perp} (g^n)^2\, d\bx^\perp = \int_{\Omega_h^\perp} (g(\bm x^\perp, t^n))^2\, d\bx^\perp,&\quad n=0,1,2,....
  \end{aligned}
\end{equation}
Surface integrals over $\Gamma_s^n = \Gamma_s(t^n)$ are interpreted analogously:
\begin{equation}
  \label{eq:weakform:surface_integral_notation_discrete}
  \int_{\Gamma_s^n} f\, ds = \int_{\Gamma_s(t^n)} f(\bm x, t^n)\, ds,\quad n=0,1,2,....
\end{equation}
}

\subsection{\reviewerTwo{Continuity regularization of the free-surface equation using a high-order edge stabilization term}}
\label{sec:discrete_problem_analysis:continuity_regularization}
Deriving a discrete counterpart to the continuous total energy estimate \eqref{eq:theorem_stability} is a crucial step 
for \reviewerTwo{attaining } \reviewerTwo{discrete $L^2$-norm stability in the total energy norm, which is 
an important building block when discussing well-posedness of the coupled problem (see Section \ref{sec:continuous_problem_analysis:well_posedness_discussion}).} 
However, \reviewerTwo{$L^2$-norm estimates of $\surf$ are not sufficient to 
ensure that $\surf$ is Lipschitz continuous, which is a requirement for well-posedness of the coupled problem as discussed in Section \ref{sec:continuous_problem_analysis:well_posedness_discussion}}. 
\reviewerTwo{To fill this gap practically, we suppress potential high-frequency oscillations by regularizing the free-surface equation 
with a high-order edge stabilization term $\edge(\surf^{n+1},w)$ \cite{Burman2004,Guermond_2013_edge}.}
\reviewerTwo{This term allows second order accuracy in space with respect to the mesh size, 
when combined with $\mathbb{P}_1$ spaces used for discretizing transport problems \cite{Guermond_2013_edge}.} 
To all the considered numerical schemes in this paper we hence add the edge stabilization term \reviewerTwo{$\edge(\surf^{n+1},w)$} defined as:
\begin{equation}
  \label{eq:discrete_problem_analysis:edge_stabilization}
\reviewerTwo{\edge(\surf^{n+1},w)} = \sum_{e \in \mathcal{E}_h^\perp} \int_e \gamma_e^n \jump{ \nabla^\perp \surf^{n+1}} \cdot \jump{ \nabla^\perp w } ds,\quad \gamma_e^n = \frac{1}{2}\,h_{K^\perp}^2 | \bu^n |,
\end{equation}
where the jump of a vector field over an interior mesh edge $e = \partial K_+ \cap \partial K_- \in \calE$ is defined as $\jump{ \bw } := \bw_+ \cdot \bn_+ - \bw_- \cdot \bn_-$ 
with $\bw_\pm = \bw|_{K_\pm}$ and $\bn_\pm$ denoting the outward-pointing unit normal vector on $\partial K_\pm$.
The mesh size $h_{K^\perp}$ is defined as the cell diameter. 

\reviewerTwo{We intentionally treat the edge stabilization term implicitly in the surface height 
to avoid inducing additional time-step restrictions. 
This does not imply large computational overheads as the term only depends on $\surf$ and is linear in $\surf$.
However, on its own, this term does not enable the $L^2$-norm stability of the free surface, and certainly not 
the unconditional $L^2$-stability with respect to $\Delta t$, which is a result we obtain later in our analysis. This is also observed in the numerical examples, 
see for instance Figure \ref{fig:experiments:greenland:numerical_solutions_diff_stabs}.

Term $\edge(\surf^{n+1}, w)$ as defined in \eqref{eq:discrete_problem_analysis:edge_stabilization} does not 
guarantee Lipschitz continuity when $\Omega^\perp \subset \mathbb{R}^2$. However, it has sufficiently improved the regularity of $\surf$ 
in the considered numerical examples. 
Precise choice of $\edge(\surf^{n+1},w)$ for the free-surface equation only recently became an active research topic \cite{bueler2024_surf_elevation_errors}. Formulations alternative to \eqref{eq:discrete_problem_analysis:edge_stabilization} are possible. 
The analysis results in the subsequent sections do not change as long 
as $\edge(\surf^{n+1},w)$ is symmetric and positive semi-definite, and thus induces a semi-norm. 
Moreover, $\edge(\surf^{n+1},w)$ should be consistent with the physics of the free-surface equation and control an appropriate norm to in the end ensure Lipschitz continuity of $\surf$ 
as discussed in Section \ref{sec:continuous_problem_analysis:well_posedness_discussion}.}
\subsection{Implicit Euler discretized problem analysis} 
\label{Sec:be_stability}
The implicit Euler discretization of \eqref{eq:stokes_free_surface_weakform}, \reviewerTwo{including symmetric edge stabilization $\edge(\surf^{n+1},w)$} \eqref{eq:discrete_problem_analysis:edge_stabilization}, 
is defined as follows. Given $\surf^n$ and $a^{n+1}$, find $\bm u^{n+1} \in \bm V_h(\Omega_h^{n+1})$, $\pi^{n+1} \in Q_h(\Omega_h^{n+1})$, 
and $\surf^{n+1} \in Z_h(\Omega_h^\perp)$ such that:
\begin{equation}
  \label{eq:backward_euler}
  \left\{
\begin{aligned}
  \lft( \frac{\surf^{n+1} - \surf^n}{\Delta t}, w \rgt)_{\Omega_h^\perp}\!\!\!\!\!\! +  (\bm u^{\perp,n+1} \cdot \nabla ^\perp \surf^{n+1} - u_z^{n+1}, w)_{\Omega_h^\perp}  + \reviewerTwo{\edge(\surf^{n+1},w)}  &= (a^{n+1},w)_{\Omega_h^\perp}, \\
  2( \sqrt{\mu}\, \bm D\bm u^{n+1}, \sqrt{\mu}\,\bm D \bm v)_{\Omega_h^{n+1}} - (\pi^{n+1}, \nabla \cdot  \bm v)_{\Omega_h^{n+1}} &= - \rho g(  \hat{\bz}, \bv)_{\Omega_h^{n+1}},\\
( \nabla \cdot  \bu^{n+1}, q)_{\Omega_h^{n+1}} &= 0,
\end{aligned}
\right.
\end{equation}
for all $\bm v \in \bm V_h(\Omega_h^{n+1})$, $q \in Q_h(\Omega_h^{n+1})$, and $w \in Z_h(\Omega_h^\perp)$. 
In the two propositions below, we show that this scheme:
\begin{enumerate}
\item is unconditionally stable with respect to $\Delta t$ in the total energy norm 
\reviewerTwo{when free-surface height is Lipschitz continuous at all times,} \reviewerOne{where 
the resulting one-step bound implies a Grönwall-type control under forcing},
\item enables domain volume conservation.
\end{enumerate}
However, the scheme requires nonlinear iterations across the Stokes/free-surface problem due to the dependence of $\Omega^{n+1}$ in the Stokes problem to $\surf^{n+1}$ in the free-surface equation. 

\begin{proposition}[Implicit Euler stability] \label{proposition:backward euler stability}
  Let $(\bm u^{n+1}, \pi^{n+1}, \surf^{n+1})$ be a solution to \eqref{eq:backward_euler} with source term $a^{n+1}=a(\bm x^\perp, t^{n+1})$. 
  \reviewerTwo{Assume that $\surf^{n+1} > b$ everywhere on $\Omega_h^\perp$.}
  Then the following stability estimate in the total energy norm holds:
  \begin{equation}\label{eq:impliciteuler_energy bound_final}
      \begin{aligned}
    \|\surf^{n+1}\|_{L^2(\Omega_h^\perp)}^2 + \frac{4}{\rho g}\, \Delta t\, \|\sqrt{\mu}\,\bm D \bm u\|_{L^2(\Omega_h^{n+1})}^2 \reviewerTwo{+ 2\, \Delta t\, \edge(\surf^{n+1}, \surf^{n+1})} &\leq \|\surf^n \|^2_{L^2(\Omega_h^\perp)}  + 2 \Delta t\, \| a^{n+1} \|_{L^2(\Omega_h^\perp)} \| \surf^{n}\|_{L^2(\Omega_h^\perp)} \\
    &\quad + 2 (\Delta t)^2 \| a^{n+1} \|_{L^2(\Omega_h^\perp)}^2.  
      \end{aligned}
    \end{equation}
\end{proposition}
\begin{proof}
\,\\\textbf{(i) Energy balance.} 
We set $w=2\surf^{n+1}$ in the free-surface equation in \eqref{eq:backward_euler} and then 
apply the identity $2\tilde a^2 - 2 \tilde a \tilde b = \tilde a^2 - \tilde b^2 + (\tilde a-\tilde b)^2$ with $\tilde a=\surf^{n+1}$ and $\tilde b=\surf^n$, to obtain: 
\begin{equation} 
  \label{eq:impliciteuler_energy bound_tmp1}
  \begin{split}
\|\surf^{n+1}\|^2_{L^2(\Omega_h^\perp)}  &= \|\surf^n \|^2_{L^2(\Omega_h^\perp)}  -\|\surf^{n+1} - \surf^n \|^2_{L^2(\Omega_h^\perp)} + 2 \Delta t (- \bm u^{\perp,n+1} \cdot \nabla^\perp \surf^{n+1} + u_z^{n+1} + a^{n+1}, \surf^{n+1})_{\Omega_h^\perp} \\ 
&\quad - 2 \Delta t\, \reviewerTwo{\edge(\surf^{n+1}, \surf^{n+1})}.
\end{split}
\end{equation}
We now focus on the term $(- \bm u^{\perp,n+1} \cdot \nabla^\perp \surf^{n+1} + u^{n+1}_z, \surf^{n+1})_{\Omega_h^\perp}$. Since all variables in this bilinear form are evaluated at the same time $t=t^{n+1}$, 
the analysis framework for the continuous case as outlined in Lemma \ref{lemma:free_surf_equation_to_surface_integral}, Lemma \ref{lemma:Stokes_to_surface_integral_ztrick}, and 
Theorem \ref{theorem:surface_height_energy_rate} applies. In particular, using \eqref{eq:lemma:surface_height_energy_rate_tmp1} from the proof of 
Theorem \ref{theorem:surface_height_energy_rate} 
we obtain: $2\,\|\sqrt{\mu}\,\bm D \bm u\|_{L^2(\Omega_h^{n+1})}^2 = -\rho g\, (- \bm u^{\perp,n+1} \cdot \nabla^\perp \surf^{n+1} + u^{n+1}_z, \surf^{n+1})_{\Omega_h^\perp}.$
Dividing the equation above by $\rho g$, substituting this result into \eqref{eq:impliciteuler_energy bound_tmp1}, and rearranging the terms, gives: 
\begin{equation}
    \label{eq:implicit_energy_balance}
    \begin{aligned}
  \|\surf^{n+1}\|_{L^2(\Omega_h^\perp)}^2   +\frac{4}{\rho g}\, \Delta t\, \|\sqrt{\mu}\,\bm D \bm u\|_{L^2(\Omega_h^{n+1})}^2 
   &= \|\surf^n \|^2_{L^2(\Omega_h^\perp)}  + 2 \Delta t\, (a^{n+1}, \surf^{n+1})_{\Omega_h^\perp} \\
   &\quad - 2\,\Delta t\, \reviewerTwo{\edge(\surf^{n+1}, \surf^{n+1})} -\|\surf^{n+1} - \surf^n \|^2_{L^2(\Omega_h^\perp)}.
    \end{aligned}
  \end{equation}
\textbf{(ii) An intermediate $L^2$-estimate of $\surf$.}
We set $w=\surf^{n+1}, \bv = \bu^{n+1}$, and $q=\pi^{n+1}$ in \eqref{eq:backward_euler}, and repeat similar steps as in the first part of this proof, to obtain: 
\begin{equation} 
  \begin{split}
    \|\surf^{n+1}\|_{L^2(\Omega_h^\perp)}^2 + \Delta t\, \reviewerTwo{\edge(\surf^{n+1}, \surf^{n+1})} &= (\surf^n, \surf^{n+1})_{\Omega_h^\perp} -\frac{2}{\rho g}\, \Delta t\, \|\sqrt{\mu}\,\bm D \bm u\|_{L^2(\Omega_h^{n+1})}^2 + \Delta t\, (a^{n+1}, \surf^{n+1})_{\Omega_h^\perp} \\ 
    &\quad - \Delta t\, \reviewerTwo{\edge(\surf^{n+1}, \surf^{n+1})}.
\end{split}
\end{equation}
Neglecting the negative terms and applying the Cauchy-Schwarz inequality, we estimate: 
\begin{equation}
\begin{aligned}
\|\surf^{n+1}\|_{L^2(\Omega_h^\perp)}^2 
 &\leq (\surf^n, \surf^{n+1})_{\Omega_h^\perp} + \Delta t\, (a^{n+1}, \surf^{n+1})_{\Omega_h^\perp} \\
 &= (\surf^n + \Delta t\, a^{n+1}, \surf^{n+1})_{\Omega_h^\perp} 
\leq  \|\surf^n + \Delta t\, a^{n+1}\|_{L^2(\Omega_h^\perp)}\, \|\surf^{n+1}\|_{L^2(\Omega_h^\perp)}.
\end{aligned}
\end{equation}

Dividing both sides by $\|\surf^{n+1}\|_{L^2(\Omega_h^\perp)}$, and applying the triangle inequality then gives:
\begin{equation} \label{eq:implicit_l2_estimate}
  \|\surf^{n+1}\|_{L^2(\Omega_h^\perp)} \leq  \|\surf^n\|_{L^2(\Omega_h^\perp)}\, + \Delta t\, \|a^{n+1}\|_{L^2(\Omega_h^\perp)}.
\end{equation}
\textbf{(iii) Total energy estimate.}
We start at \eqref{eq:implicit_energy_balance}, \reviewerTwo{move $-2 \Delta t\, \reviewerTwo{\edge(\surf^{n+1}, \surf^{n+1})}$ to the left-hand-side}, and use $2 \Delta t\, (a^{n+1}, \surf^{n+1})_{\Omega_h^\perp} \leq 2 \Delta t\, \| a^{n+1} \|_{L^2(\Omega_h^\perp)} \| \surf^{n+1}\|_{L^2(\Omega_h^\perp)}$. 
Then we have: 
\begin{equation}
  \begin{aligned}
&\|\surf^{n+1}\|_{L^2(\Omega_h^\perp)}^2 + \frac{4}{\rho g}\, \Delta t\, \|\sqrt{\mu}\,\bm D \bm u\|_{L^2(\Omega_h^{n+1})}^2 +  \reviewerTwo{ 2\, \Delta t\edge(\surf^{n+1}, \surf^{n+1})}   \\
&\leq \|\surf^n \|^2_{L^2(\Omega_h^\perp)}  + 2 \Delta t\, \| a^{n+1} \|_{L^2(\Omega_h^\perp)} \| \surf^{n+1}\|_{L^2(\Omega_h^\perp)} \\
  &\leq \|\surf^n \|^2_{L^2(\Omega_h^\perp)}  + 2 \Delta t\, \| a^{n+1} \|_{L^2(\Omega_h^\perp)} \| \surf^{n}\|_{L^2(\Omega_h^\perp)} + 2 (\Delta t)^2 \| a^{n+1} \|_{L^2(\Omega_h^\perp)}^2,  
  \end{aligned}
\end{equation}
where we used \eqref{eq:implicit_l2_estimate} in the final step.
\end{proof}

\begin{proposition}[Implicit Euler conservation] 
  \label{proposition:backward_euler_domain_conservation}
   The solution to \eqref{eq:backward_euler} with source term $a^{n+1}=a(\bm x^\perp, t^{n+1})$ 
   \reviewerTwo{and assumption $\surf^{n+1} > b$ everywhere on $\Omega_h^\perp$,} conserves domain volume: $|\Omega_h^{n+1}| = |\Omega_h^n| + \Delta t\, \int_{\Omega_h^\perp} a^{n+1}\, d \bx^\perp.$
\end{proposition}
\begin{proof}
We set $w=1$ in the free-surface equation part of \eqref{eq:backward_euler} to obtain:
\begin{equation}
  \label{eq:impliciteuler_volume_conservation_tmp1}
(\surf^{n+1}, 1)_{\Omega_h^\perp} = (\surf^n, 1)_{\Omega_h^\perp} + \Delta t\, (- \bm u^{\perp,n+1} \cdot \nabla^\perp \surf^{n+1} + u^{n+1}_z + a^{n+1}, 1)_{\Omega_h^\perp},
\end{equation}
\reviewerOne{where we have neglected the edge stabilization term $\edge(\surf^{n+1}, 1)$ as it vanishes when $w=1$.}
Since all terms in the last integrand are evaluated at time $t = t^{n+1}$, the analysis in Section \ref{Sec:volume_conservation} applies. 
Following the derivation of the continuous volume conservation in \eqref{eq:continuous_volume_conservation_tmp2}, we have that: 
\begin{equation}
  \label{eq:impliciteuler_volume_conservation_tmp2}
(- \bm u^{\perp,n+1} \cdot \nabla^\perp \surf^{n+1} + u^{n+1}_z, 1)_{\Omega_h^\perp} = \int_{ \partial\Omega_h^{n+1}} (\bm u \cdot \bm n)\, ds = \int_{\Omega_h^{n+1}} \nabla \cdot \bm u^{n+1}\, d \bx = 0,
\end{equation}
where the final equality follows from the incompressibility condition \eqref{eq:backward_euler} with $q=1 \in Q_h(\Omega_h^{n+1})$. Combining \eqref{eq:impliciteuler_volume_conservation_tmp1} and \eqref{eq:impliciteuler_volume_conservation_tmp2}, and using the definition of domain volume via the height function, we conclude
$|\Omega_h^{n+1}| - |\Omega_h^n| = (\surf^{n+1} -b,1)_{\Omega_h^\perp} - (\surf^{n} -b,1)_{\Omega_h^\perp} = (\surf^{n+1} ,1)_{\Omega_h^\perp} - (\surf^{n},1)_{\Omega_h^\perp} =   \Delta t\, \int_{\Omega_h^\perp} a^{n+1}\, d \bx^\perp,$
which completes the proof.
\end{proof}

\subsection{\reviewerTwo{Explicit Euler with implicit surface regularization}}
\label{sec:discrete_analysis:explicit_euler}
We now consider the explicit Euler discretization of \eqref{eq:stokes_free_surface_weakform} with \reviewerTwo{implicit symmetric edge stabilization $\edge(\surf^{n+1},w)$} \eqref{eq:discrete_problem_analysis:edge_stabilization}. 
Given $\surf^n$ and $a^n$, find $\bm u^n \in \bm V_h(\Omega_h^n)$, $\pi^n \in Q_h(\Omega_h^n)$, and $\surf^{n+1} \in Z_h(\Omega_h^\perp)$ such that:
\begin{equation}
  \label{eq:forward_euler}
  \left\{
\begin{aligned}
  \lft( \frac{\surf^{n+1} - \surf^n}{\Delta t}, w \rgt)_{\Omega_h^\perp} +  \lft( \bm u^{\perp,n} \cdot \nabla ^\perp \surf^{n} - u_z^{n}, w \rgt)_{\Omega_h^\perp} + \reviewerTwo{\edge(\surf^{n+1},w)}  &= (a^n,w)_{\Omega_h^\perp}, \\
  2 \lft( \sqrt{\mu}\,\bm D\bm u^{n}, \sqrt{\mu}\,\bm D \bm v \rgt)_{\Omega_h^{n}} - (\pi^{n}, \nabla \cdot  \bm v)_{\Omega_h^{n}} &= - \rho g  ( \hat{\bz}, \bv)_{\Omega_h^{n}},\\
( \nabla \cdot  \bu^{n}, q)_{\Omega_h^{n}} &= 0,
\end{aligned}
\right.
\end{equation}
for all $\bm v \in \bm V_h(\Omega_h^{n})$, $q \in Q_h(\Omega_h^{n})$, and $w \in Z_h(\Omega_h^\perp)$. 

This is a numerical scheme central to this paper. The scheme does not require 
nonlinear iterations across the Stokes/free-surface problems 
\reviewerTwo{which are the main computational burden when using implicit methods such as \eqref{eq:backward_euler}. 
Scheme \eqref{eq:forward_euler}, however, requires linear implicit solves at the free-surface equation level due to: (i) non-diagonal mass matrix, (ii) the 
implicit surface regularization. These implicit solves are computationally less expensive compared to the nonlinear iterations across the Stokes/free-surface problems.} 
\reviewerTwo{In Proposition \ref{proposition:stability_forward_euler} and Proposition \ref{proposition:forward_euler_domain_volume_cons}, stated below,} we show that the scheme \eqref{eq:forward_euler}:
\begin{enumerate}
\item is unstable in the total energy norm for any $\Delta t$, with spurious terms $(\Delta t)^2\,\int_{\Gamma_s^n} \omega\, (\bm u \cdot \bm n)^2\, ds$ and  $2\,(\Delta t)^2\, \int_{\Gamma_s^n} (\bm u \cdot \bm n)\, a^n\, ds$ 
that \reviewerOne{prevent a discrete analogue of the continuous total energy estimate \eqref{eq:theorem_stability}},
\item enables domain volume conservation.
\end{enumerate}
\reviewerTwo{While the two spurious terms are solution dependent, they are not a part of the total energy norm as shown in Proposition \ref{proposition:stability_forward_euler}. Moreover, one of them is positive, whereas the sign of the other 
is not fixed. Together, the terms cause the total energy to grow in time and 
thus leads to an unstable solution as also demonstrated later using numerical examples.}

\reviewerOne{We note that this type of instability with respect to $\Delta t$, which arises from the free-surface equation, 
carries similarities with the well known instability of the linear advection equation discretized by central differences in space and explicit Euler 
discretization in time. The instability is only related to treating the free-surface equation explicitly and is not 
related to for instance lack of regularity of the solution, as the edge stabilization term is treated implicitly and thus does not induce additional time-step restrictions.}

\begin{proposition}[Explicit Euler stability estimate] 
  \label{proposition:stability_forward_euler}
  Let $(\bm u^n, \pi^n, \surf^{n+1})$ be a solution to \eqref{eq:forward_euler} with source term $a^n=a(\bm x^\perp, t^n)$. 
  \reviewerTwo{Assume $\surf^{n+1} > b$ everywhere on $\Omega_h^\perp$.} Then the following stability estimate in the total energy norm holds:
  \begin{equation} \label{eq:explicit_stability_final}
    \begin{aligned}
      &\|\surf^{n+1}\|_{L^2(\Omega_h^\perp)}^2 + \frac{4}{\rho g}\, \Delta t\, \|\sqrt{\mu}\, \bm D \bm u\|^2_{L^2(\Omega_h^n)} \reviewerTwo{+ 2\, \Delta t\, \edge(\surf^{n+1}, \surf^{n+1})}\\
      &\leq  \|\surf^{n}\|_{L^2(\Omega_h^\perp)}^2 + 2\, \Delta t(a^n, \surf^n)_{\Omega^\perp_h} +(\Delta t)^2 \|a^n\|^2_{L^2(\Omega^\perp_h)} 
       +  (\Delta t)^2\,\int_{\Gamma_s^n} \omega\, (\bm u \cdot \bm n)^2\, ds + 2\,(\Delta t)^2\, \int_{\Gamma_s^n} (\bm u \cdot \bm n)\, a^n\, ds,
    \end{aligned}
  \end{equation}
  \reviewerTwo{with $\omega = \| ( -\nabla^\perp \surf, 1) \|_2$}.
 \end{proposition}
  \begin{proof}
\,\\
\textbf{(i) Decomposition of the free-surface equation into explicit and implicit parts.}
In the first equation of \eqref{eq:forward_euler}, we \reviewerTwo{substitute $\surf^{n+1} = \surf_E^{n+1} + \surf_I^{n+1}$ and 
rearrange the terms to arrive at: 
$$\lft( \frac{\surf_E^{n+1} - \surf^n}{\Delta t}, w \rgt)_{\Omega_h^\perp} +  (\bm u^{\perp,n} \cdot \nabla ^\perp \surf^{n} - u_z^{n}, w)_{\Omega_h^\perp}  = (a^n, w)_{\Omega_h^\perp} - \lft( \frac{\surf_I^{n+1}}{\Delta t}, w \rgt)_{\Omega_h^\perp} - \reviewerTwo{\edge(\surf^{n+1},w)}.$$ 
This equation is a sum of the following two equations:
    \begin{equation} \label{eq:hom-update}
      \lft( \frac{\surf_E^{n+1} - \surf^n}{\Delta t}, w \rgt)_{\Omega_h^\perp} +  (\bm u^{\perp,n} \cdot \nabla ^\perp \surf^{n} - u_z^{n}, w)_{\Omega_h^\perp}  = (a^n, w)_{\Omega_h^\perp},
    \end{equation}
    \begin{equation} \label{eq:edge-update}
      \lft( \frac{\surf_I^{n+1}}{\Delta t}, w \rgt)_{\Omega_h^\perp}   =  -\reviewerTwo{\edge(\surf^{n+1},w)},
    \end{equation}
which is a consistent decomposition into two separate forms containing explicit terms and implicit terms, used to estimate the energy contributions 
of the explicit and implicit parts of the free-surface equation separately.
}
In \eqref{eq:hom-update}, we now set $w = 2\surf^n$ and use the identity $2 \tilde a \tilde b - 2 \tilde b^2 = \tilde a^2 - \tilde b^2 - (\tilde a-\tilde b)^2$ with $\tilde a=\surf_E^{n+1}$, $\tilde b=\surf^n$, to obtain: 
  \begin{equation}
    \label{eq:explicit_stability_tmp1}
    \begin{aligned}
  \|\surf_E^{n+1}\|_{L^2(\Omega_h^\perp)}^2 &=  \|\surf^{n}\|_{L^2(\Omega_h^\perp)}^2 + 2 \Delta t (- \bm u^{\perp,n} \cdot \nabla ^\perp \surf^{n} + u_z^n, \surf^n)_{\Omega_h^\perp} + 2 \Delta t (a^n, \surf^n)_{\Omega_h^\perp}   + \|\surf_E^{n+1} - \surf^n\|^2_{L^2(\Omega_h^\perp)}.
    \end{aligned}
  \end{equation}
 \textbf{(ii) Energy balance.}
  We use Lemma \ref{lemma:free_surf_equation_to_surface_integral} to write term $(- \bm u^{\perp,n} \cdot \nabla ^\perp \surf^{n} + u_z^n, \surf^n)_{\Omega_h^\perp}$ 
  as a surface integral, 
  followed by using Lemma \ref{lemma:Stokes_to_surface_integral_ztrick} to write the resulting surface integral to the scaled Stokes gravity term. This gives:
  $2\, \Delta t (- \bm u^{\perp,n} \cdot \nabla^\perp \surf^n + u_z^n,\surf^n)_{\Omega_h^\perp} = 2\, \Delta t\, \int_{\Gamma_s^n} (\bm u \cdot \bm n)\, z\, ds
  =  2 \Delta t\, ( \hat{\bz}, \bm u)_{\Omega_h^n}$. 
  Next, we use the Stokes problem \eqref{eq:forward_euler} formulation with $\bm v = \bm u^n$ and $q = \pi^n$ to relate the scaled gravity term to the strain tensor norm as:
  \begin{equation}
    \label{eq:explicit_stability_tmp1_1}
    2\, \Delta t (- \bm u^{\perp,n} \cdot \nabla^\perp \surf^n + u_z^n,\surf^n)_{\Omega_h^\perp} =  2 \Delta t\, (\hat{\bz}, \bm u)_{\Omega_h^n} = -\frac{4}{\rho g}\, \Delta t\, \|\sqrt{\mu}\, \bm D \bm u\|^2_{L^2(\Omega_h^n)}.
  \end{equation}
  Inserting this relation into \eqref{eq:explicit_stability_tmp1} and then rearranging, gives the energy balance:
  \begin{equation}
    \label{eq:explicit_stability_tmp2}
    \begin{aligned}
  \| \surf_E^{n+1} \|_{L^2(\Omega_h^\perp)}^2 + \frac{4}{\rho g}\, \Delta t\, \|\sqrt{\mu}\, \bm D \bm u\|^2_{L^2(\Omega_h^n)} &=  \|\surf^{n}\|_{L^2(\Omega_h^\perp)}^2  + \|\surf_E^{n+1} - \surf^n\|^2_{L^2(\Omega_h^\perp)} +  2 \Delta t (a^n, \surf^n)_{\Omega_h^\perp}.
    \end{aligned}
  \end{equation}
\textbf{(iii) Spurious term estimate.} 
  We estimate the spurious term $\|\surf_E^{n+1} - \surf^n\|^2_{L^2(\Omega_h^\perp)}$ in \eqref{eq:explicit_stability_tmp2} by setting $w = \surf_E^{n+1} - \surf^n$ and then estimating: 
    $\| \surf_E^{n+1} - \surf^n \|_{L^2(\Omega_h^\bot)}^2 = \Delta t ( -\bm u^{\perp,n} \cdot \nabla ^\perp \surf^{n} + u_z^{n} + a^n, \surf_E^{n+1} - \surf^n )_{\Omega_h^\bot} 
    \leq \Delta t\, \| -\bm u^{\perp,n} \cdot \nabla ^\perp \surf^{n} + u_z^{n} + a^n\|_{L^2(\Omega_h^\bot)} \, \| \surf_E^{n+1} - \surf^n \|_{L^2(\Omega_h^\bot)}.$
  Dividing both sides by $\| \surf_E^{n+1} - \surf^n \|_{L^2(\Omega_h^\perp)}$ and squaring the resulting inequality gives:
  \begin{equation} \label{eq:spurious}
    \begin{aligned}
  \| \surf_E^{n+1} - \surf^n \|^2_{L^2(\Omega_h^\perp)} &\leq (\Delta t)^2 \, \| -\bm u^{\perp,n} \cdot \nabla ^\perp \surf^{n} + u_z^{n} + a^n \|^2_{L^2(\Omega_h^\bot)} \\
  &= (\Delta t)^2 \, \underbrace{\| -\bm u^{\perp,n} \cdot \nabla ^\perp \surf^{n} + u_z^{n} \|^2_{L^2(\Omega_h^\perp)}}_{\text{\textbf{Term A}}}  \\
  &\quad + 2(\Delta t)^2\, \underbrace{(a^n, -\bm u^{\perp,n} \cdot \nabla ^\perp \surf^{n} + u_z^{n})_{\Omega_h^\perp}}_{\text{\textbf{Term B}}} + (\Delta t)^2\, \|a^n \|^2_{L^2(\Omega_h^\bot)}.
    \end{aligned}
  \end{equation}
  We now rewrite Term A in \eqref{eq:spurious} by defining an integral scaling: 
  $\omega = \| ( -\nabla^\perp \surf^n, 1 ) \|_2 = \frac{1}{n_z} > 0$, 
  where $n_z$ denotes the vertical component of an outward unit normal vector $\bn$ on $\Gamma_s^n$.
 Inserting this relation into \eqref{eq:spurious}, we have: 
  \begin{equation} \label{eq:term_A}
    \begin{aligned}
  \| -\bm u^{\perp,n} \cdot \nabla ^\perp \surf^{n} + u_z^{n} \|^2_{L^2(\Omega_h^\bot)}  &= \int_{\Omega_h^\perp} \hspace{-0.1cm}\frac{\omega}{\omega}\, (-\bm u^{\perp,n} \cdot \nabla ^\perp \surf^{n} + u_z^{n})\, (-\bm u^{\perp,n} \cdot \nabla ^\perp \surf^{n} + u_z^{n})\, d \bx^\perp \\
  &= \int_{\Gamma_s^n} \hspace{-0.1cm}(\bm u \cdot \bm n)\,(-\bm u^{\perp,n} \cdot \nabla ^\perp \surf^{n} + u_z^{n})\, ds = \int_{\Gamma_s^n} \hspace{-0.1cm}\omega\,(\bm u \cdot \bm n)^2 ds,
    \end{aligned}
  \end{equation}
  where we in the second step used: $\frac{1}{\omega}\, (-\bm u^{\perp,n} \cdot \nabla ^\perp \surf^{n} + u_z^{n}) = \bm u^n \cdot \bm n$ and $\omega\, d \bx^\perp = ds$ (as in the proof of Lemma \ref{lemma:free_surf_equation_to_surface_integral}). 
  In the final step, we used $-\bm u^{\perp,n} \cdot \nabla ^\perp \surf^{n} + u_z^{n} = \omega\, (\bm u \cdot \bm n)$. 
  We proceed by rewriting Term B in \eqref{eq:spurious} as follows:
  \begin{equation}
    \label{eq:term_B}
  (a^n, -\bm u^{\perp,n} \cdot \nabla ^\perp \surf^{n} + u_z^{n})_{\Omega_h^\perp} = \int_{\Omega_h^\perp} \frac{\omega}{\omega} a^n\, (-\bm u^{\perp,n} \cdot \nabla ^\perp \surf^{n} + u_z^{n})\, d \bx^\perp = \int_{\Gamma_s^n} (\bm u \cdot \bm n)\, a^n\, ds.\\
  \end{equation}
  Altogether, inserting \eqref{eq:term_A} and \eqref{eq:term_B} into \eqref{eq:spurious}, gives the final bound of the spurious term:
\begin{equation} \label{eq:spurious_estimate}
  \| \surf_E^{n+1} - \surf^n \|_{L^2(\Omega_h^\bot)}^2 \leq (\Delta t)^2\, \int_{\Gamma_s^n} \omega\, (\bm u \cdot \bm n)^2\, ds + 2\, (\Delta t)^2\, \int_{\Gamma_s^n} (\bm u \cdot \bm n)\, a^n\, ds + (\Delta t)^2\, \|a^n \|^2_{L^2(\Omega_h^\bot)}.
\end{equation}
\textbf{(iv) Finalization.} 
We now make steps towards estimating $\|\surf^{n+1}\|_{L^2(\Omega^\perp_h)}^2$ by writing 
$\|\surf^{n+1}\|_{L^2(\Omega^\perp_h)}^2 = \| \surf_E^{n+1} + \surf_I^{n+1} \|^2_{L^2(\Omega_h^\bot)} 
= \| \surf_E^{n+1} \|^2_{L^2(\Omega_h^\bot)} + \| \surf_I^{n+1}  \|^2_{L^2(\Omega_h^\bot)} + 2(\surf_E^{n+1}, \surf_I^{n+1})_{\Omega_h^\perp}$. 
We continue by testing 
\eqref{eq:edge-update} against $w=\surf_E^{n+1} + \surf_I^{n+1}$ to first 
expand the term \reviewerTwo{$2(\surf_E^{n+1}, \surf_I^{n+1})_{\Omega_h^\perp} = -2 \|\surf_I^{n+1}\|^2 - 2\, \Delta t\, \edge(\surf^{n+1}, \surf^{n+1})$ to arrive at:} 
\begin{equation}
  \begin{aligned}
  \label{eq:explicit_euler_tmp_explicit_implicit_tmp1}
\|\surf^{n+1}\|_{L^2(\Omega^\perp_h)}^2  &= \| \surf_E^{n+1} \|^2_{L^2(\Omega_h^\bot)} + \|\surf_I^{n+1} \|^2_{L^2(\Omega_h^\perp)} \reviewerTwo{- 2\|\surf_I^{n+1} \|^2_{L^2(\Omega_h^\perp)} - 2\,\Delta t\, \edge(\surf^{n+1}, \surf^{n+1})} \\
&\leq \reviewerTwo{\| \surf_E^{n+1} \|^2_{L^2(\Omega_h^\bot)} - 2\,\Delta t\, \edge(\surf^{n+1}, \surf^{n+1})}.
  \end{aligned}
\end{equation}
After rearranging, it follows that $\|\surf^{n+1}\|_{L^2(\Omega^\perp_h)}^2 \reviewerTwo{+ 2\,\Delta t\, \edge(\surf^{n+1}, \surf^{n+1})  \leq} \| \surf_E^{n+1} \|^2_{L^2(\Omega_h^\bot)}$.
Finally, inserting \eqref{eq:explicit_stability_tmp2} into the latter relation, followed by 
inserting \eqref{eq:spurious_estimate} into the resulting estimate, proves this proposition.
  \end{proof}

\begin{proposition}[Explicit Euler domain volume conservation] \label{proposition:forward_euler_domain_volume_cons}
  Solution to the explicit Euler scheme \eqref{eq:forward_euler} with source term $a^n=a(\bm x^\perp, t^n)$ 
  \reviewerTwo{and assumption $\surf^{n+1} > b$ everywhere on $\Omega_h^\perp$,} 
  conserves domain volume: $|\Omega_h^{n+1}| - |\Omega_h^n| = \Delta t\, \int_{\Omega^\perp_h} a^{n} \ud \bx^\perp$.
\end{proposition}

\begin{proof}
Applying all steps of the proof of Proposition \ref{proposition:backward_euler_domain_conservation} to \eqref{eq:forward_euler}, proves this proposition.
\end{proof}

  \subsubsection{\reviewerTwo{Weakly stable time step size for explicit Euler with implicit surface regularization}}
  The stability estimate \eqref{eq:explicit_stability_final} implies that the corresponding scheme \eqref{eq:forward_euler} is unstable 
  in the total energy norm for any $\Delta t$. 
  In the steps that follow, we determine $\Delta t$ that implies stability in a weaker sense, that is, stability without having a full control over the total energy norm. 
 This is disadvantageous as it implies a weaker coupling between the Stokes/free-surface equations and gives rise to increased energy errors. 
 \reviewerTwo{In the steps that follow, we explore and formalize this regime as many practitioners use \eqref{eq:forward_euler} and then decrease $\Delta t$ to obtain a somewhat stable solution, yet 
 the solution is not stable in the total energy norm, which is undesirable from a physical point of view.}

 We first write 
  $\frac{4}{\rho g}\, \Delta t\, \|\sqrt{\mu}\, \bm D \bm u\|^2_{L^2(\Omega_h^n)} = (1 - \varepsilon)\,\frac{4}{\rho g}\, \Delta t\,\|\sqrt{\mu}\, \bm D \bm u\|^2_{L^2(\Omega_h^n)} + \varepsilon\, \frac{4}{\rho g}\, \Delta t\, \|\sqrt{\mu}\, \bm D \bm u\|^2_{L^2(\Omega_h^n)}$, 
  where $0 < \varepsilon < 1$. 
  We insert that relation into \eqref{eq:explicit_stability_final} and rearrange the terms to arrive at a weaker stability estimate:
  \begin{equation}
    \begin{aligned}
      \label{eq:explicit_stability_final_weakerstab_tmp0}
      &\|\surf^{n+1}\|_{L^2(\Omega_h^\perp)}^2 + (1-\varepsilon)\, \frac{4}{\rho g}\, \Delta t\, \|\sqrt{\mu}\, \bm D \bm u\|^2_{L^2(\Omega_h^n)} \reviewerTwo{+ 2\, \Delta t\, \edge(\surf^{n+1}, \surf^{n+1})} \\
      &\leq  \|\surf^{n}\|_{L^2(\Omega_h^\perp)}^2 + 2\, \Delta t(a^n, \surf^n)_{\Omega^\perp_h} +(\Delta t)^2 \|a^n\|^2_{L^2(\Omega^\perp_h)} \\
      &\quad +  (\Delta t)^2\,\int_{\Gamma_s^n} \omega\, (\bm u \cdot \bm n)^2\, ds + 2\,(\Delta t)^2\, \int_{\Gamma_s^n} (\bm u \cdot \bm n)\, a^n\, ds - \varepsilon\, \frac{4}{\rho g}\, \Delta t\, \|\sqrt{\mu}\, \bm D \bm u\|^2_{L^2(\Omega_h^n)}.
    \end{aligned}
  \end{equation}
For weaker stability we require that:
$$(\Delta t)^2\,\int_{\Gamma_s^n} \omega\, (\bm u \cdot \bm n)^2\, ds + 2\,(\Delta t)^2\, \int_{\Gamma_s^n} (\bm u \cdot \bm n)\, a^n\, ds - \varepsilon\, \frac{4}{\rho g}\, \Delta t\, \|\sqrt{\mu}\, \bm D \bm u\|^2_{L^2(\Omega_h^n)} \leq 0.$$ 
Solving that inequality for $\Delta t$ implies a time step restriction that gives stability in a weaker sense:
\begin{equation} \label{eq:explicit_stability_final_weakerstab_tmp1}
  \begin{aligned}
\Delta t^n \leq \varepsilon\,\frac{ \frac{4}{\rho g}\,\|\sqrt{\mu}\, \bm D \bm u\|^2_{L^2(\Omega_h^n)}}{\int_{\Gamma_s^n} \omega\, (\bm u \cdot \bm n)^2\, ds + 2\, \int_{\Gamma_s^n} (\bm u \cdot \bm n)\, a^n\, ds},\quad 0 < \varepsilon < 1,\quad n=0,1,2,...
  \end{aligned}
\end{equation}
\igor{Using such $\Delta t$,} the norm $(1-\varepsilon)\, \frac{4}{\rho g}\, \Delta t\, \|\sqrt{\mu}\, \bm D \bm u\|^2_{L^2(\Omega_h^n)}$ in \eqref{eq:explicit_stability_final_weakerstab_tmp0} 
is now controlled only up to the $\varepsilon$ parameter. 
\igor{We use $\Delta t$ computed according to \eqref{eq:explicit_stability_final_weakerstab_tmp1} only for drawing comparisons across numerical schemes in the numerical experiments section. We do not 
recommend using this approach in practice} \igorr{as it leads to a weak coupling between the Stokes/free-surface equations and increased energy errors.}
Instead, in Section \ref{sec:discrete_analysis:explicit_euler_new_stabilization}, we give a stabilization term that (i) allows a full total energy norm control without $\varepsilon$-dependence, and (ii) enables stability for any choice of $\Delta t$.

\subsection{\reviewerTwo{Explicit Euler with coupled-energy stabilization and implicit surface regularization}}
\label{sec:discrete_analysis:explicit_euler_new_stabilization}
Motivated by the stability estimate \eqref{eq:explicit_stability_final} for the explicit Euler method \eqref{eq:forward_euler}, we introduce the following stabilized method.
Given $\surf^n$ and $a^n$, find $\surf^{n+1} \in Z_h(\Omega_h^\perp)$, $\bm u^n \in \bm V_h(\Omega_h^n)$, $\pi^n \in Q_h(\Omega_h^n)$ such that: 
\begin{equation}
  \label{eq:forward_euler_stabilized}
  \left\{
\begin{aligned}
  \lft( \frac{\surf^{n+1} - \surf^n}{\Delta t}, w \rgt)_{\Omega_h^{\perp}} +  (\bm u^{\perp,n} \cdot \nabla ^\perp \surf^{n} - u_z^{n}, w)_{\Omega_h^{\perp}}  + \reviewerTwo{\edge(\surf^{n+1},w)} &= (a^{n},w)_{\Omega_h^{\perp}}, \\
  \!2(\sqrt{\mu}\, \bm D\bm u^{n}, \sqrt{\mu}\, \bm D \bm v)_{\Omega_h^{n}} + \frac{ \rho g \Delta t }{2} \, \int_{\Gamma_s^n} \!\!\!\!\omega\, (\bm u \cdot \bm n)\, (\bm v \cdot \bm n)\, ds  - (\pi^{n}, \nabla \cdot  \bm v)_{\Omega_h^{n}} &= - \rho g(  \hat{\bz}, \bv)_{\Omega_h^{n}} \reviewerTwo{-  \rho g \, \Delta t\, \int_{\Gamma_s^n} \!\!\!\!(\bm v \cdot \bm n)\,a\, ds},\\
( \nabla \cdot  \bu^{n}, q)_{\Omega_h^{n}} &= 0,
\end{aligned}
\right.
\end{equation}
for all $w \in Z_h(\Omega_h^\perp)$, $\bm v \in \bm V_h(\Omega_h^{n})$, and $q \in Q_h(\Omega_h^{n})$, 
\reviewerTwo{where $\edge(\surf^{n+1},w)$ is the symmetric edge stabilization term defined in \eqref{eq:discrete_problem_analysis:edge_stabilization} 
and $\omega = \| ( -\nabla^\perp \surf^n, 1) \|_2$.} 
\reviewerTwo{This stabilized scheme:
\begin{enumerate}
\item does not require nonlinear iterations across the Stokes/free-surface problem, and
\item requires only linear implicit solves at the free-surface equation level (due to standard non-diagonal FEM mass matrix and added 
implicit surface regularization) which is the same as in the unstabilized explicit Euler method \eqref{eq:forward_euler} case.
\end{enumerate} 
}
The only difference between 
\reviewerTwo{schemes \eqref{eq:forward_euler} and \eqref{eq:forward_euler_stabilized} are the two added stabilization terms}
$\frac{ \rho g \Delta t }{2} \, \int_{\Gamma_s^n} \omega\, (\bm u \cdot \bm n)\, (\bm v \cdot \bm n)\, ds$, and $\rho g \, \Delta t\, \int_{\Gamma_s^n} a\, (\bm v \cdot \bm n)\, ds$. 
The first term is, up to the scaling, the same as introduced in \cite{rose_buffet_heister_u_dot_n}, 
where it was used for time step size estimation. 
The scaling that we propose 
is motivated by our analysis in Proposition \ref{proposition:stability_forward_euler} without simplifications 
of the original Stokes/free-surface problem formulation. 
The added stabilization terms can be seen as of FSSA-type \cite{Kaus_fssa}, 
however, compared to the FSSA stabilization terms in \eqref{eq:semi_backward_euler_fssa_stabilized}, our stabilization terms are derived such that they precisely cancel the spurious term in the 
total energy estimate outlined in Proposition \ref{proposition:stability_forward_euler}. 
\reviewerOne{The relation between our new stabilization used in \eqref{eq:forward_euler_stabilized} 
and FSSA \cite{Kaus_fssa} is: 
\begin{align}
 \frac{ \rho g \Delta t }{2} \, \int_{\Gamma_s^n} \!\!\!\!\omega\, (\bm u \cdot \bm n)\, (\bm v \cdot \bm n)\, ds 
 &= \frac{ \rho g \Delta t }{2} \, \int_{\Gamma_s^n} \!\!\!\! \|(-\nabla^\perp \surf^n,1)\|_2  \, (\bm u \cdot \bm n)\, (\bm v \cdot (-\nabla^\perp \surf^n,1)/\|(-\nabla^\perp \surf^n,1)\|_2 )\, ds \nonumber \\  
 & = \frac{ \rho g \Delta t }{2} \, \int_{\Gamma_s^n} \!\!\!\!   \, (\bm u \cdot \bm n)\, (\bm v \cdot \hat{\bz} )\, ds +  \frac{ \rho g \Delta t }{2} \, \int_{\Gamma_s^n} \!\!\!\!   \, (\bm u \cdot \bm n)\, (\bm v \cdot  (-\nabla^\perp \surf^n,0) )\, ds,
\end{align}
where $\frac{ \rho g \Delta t }{2} \, \int_{\Gamma_s^n} \!   \, (\bm u \cdot \bm n)\, (\bm v \cdot  (-\nabla^\perp \surf^n,0) )\, ds$ in the end 
is the correction to the FSSA stabilization term $\frac{ \rho g \Delta t }{2} \, \int_{\Gamma_s^n} \!   \, (\bm u \cdot \bm n)\, (\bm v \cdot \hat{\bz} )\, ds$, required to obtain the precise cancellation of the spurious term in the energy estimate.
Our new stabilization terms do not break the symmetry of the Stokes problem, which is advantageous when using iterative linear solvers.}
\reviewerOne{The key idea behind this stabilization is to hide stabilization terms in the Stokes problem, and then 
seek their cancellation with the spurious term in the energy estimate. 
This is in contrast to the approach of adding explicit stabilization terms directly to the free-surface equation, which does 
not allow such cancellation and leads to 
a CFL condition.} 

\reviewerTwo{In Proposition \ref{proposition:stability_forward_euler_stabilized} and Proposition \ref{proposition:forward_euler_stabilized_domain_volume_cons}, stated below,} we show that the scheme \eqref{eq:forward_euler_stabilized}:
\begin{enumerate}
\item is unconditionally stable with respect to $\Delta t$ in the total energy norm \reviewerTwo{with implicit surface regularization term,} \reviewerOne{where 
the resulting one-step bound implies a Grönwall-type control under forcing,}
\item enables domain volume conservation.
\end{enumerate}
\reviewerOne{As shown in the proposition, the unconditional stability is a consequence of the special structure given by the coupling of 
the Stokes and free-surface equations, and the particular choice of the stabilization terms.}
\reviewerTwo{Error estimates taking into account the added stabilization terms are out of the scope of this paper. 
Solution errors from a numerical perspective are provided in the experiments section.}
\begin{proposition}[Stabilized explicit Euler stability estimate] \label{proposition:stability_forward_euler_stabilized}
  Let $(\bm u^n, \pi^n, \surf^{n+1})$ be a solution to the \reviewerTwo{explicit Euler scheme with coupled-energy stabilization and implicit surface regularization} \eqref{eq:forward_euler_stabilized} 
  with  source term $a^n=a(\bm x^\perp, t^n)$. 
  \reviewerTwo{Assume $\surf^{n+1} > b$ everywhere on $\Omega_h^\perp$.} Then the following stability estimate in the total energy norm holds:
  \begin{equation} \label{eq:stabilized_stability_estimate}
    \begin{split}
      \| \surf^{n+1} \|_{L^2(\Omega_h^\perp)}^2 + \frac{4\Delta t}{\rho g}\, \, \|\sqrt{\mu}\, \bm D \bm u \|^2_{L^2(\Omega_h^n)} \reviewerTwo{+ 2\, \Delta t\, \edge(\surf^{n+1}, \surf^{n+1})}
      \leq  \|\surf^{n}\|_{L^2(\Omega_h^\perp)}^2 + 2\, \Delta t(a^n, \surf^n)_{\Omega^\perp_h} + (\Delta t)^2 \|a^n\|^2_{L^2(\Omega^\perp_h)}      .
    \end{split}
  \end{equation}
\end{proposition}
  
  \begin{proof}
We reuse the proof of Proposition \ref{proposition:stability_forward_euler}. Taking 
into account the stabilization terms in \eqref{eq:forward_euler_stabilized} we rewrite \eqref{eq:explicit_stability_tmp1_1} to:
  $2 \Delta t (-\bm u^{\perp,n} \cdot \nabla ^\perp \surf^{n} + u_z^{n}, \surf^n)_{\Omega_h^{\perp}} = 2\, \Delta t\, ( \hat{\bz}, \bm u)_{\Omega_h^n} 
  = -\frac{4  \Delta t}{\rho g}\,\, \|\sqrt{\mu}\, \bm D \bm u\|^2_{L^2(\Omega_h^n)} 
  - (\Delta t)^2  \, \int_{\Gamma_s^n} \omega\, (\bm u \cdot \bm n)^2\, ds - 2\,(\Delta t)^2 \, \int_{\Gamma_s^n} (\bm u \cdot \bm n)\,a^n\, ds. $ 
We now use this relation to reiterate the remaining steps in the proof of Proposition \ref{proposition:stability_forward_euler}  (from \eqref{eq:explicit_stability_tmp1_1} on) 
and then obtain the following estimate:
\begin{equation} \label{eq:messy_energy_balance2}
  \begin{split}
    &\| \surf^{n+1} \|^2_{L^2(\Omega_h^\bot)} + \frac{4 \Delta t}{\rho g }\, \, \|\sqrt{\mu}\, \bm D \bm u\|^2_{L^2(\Omega_h^n)} \reviewerTwo{+ 2\, \Delta t\, \edge(\surf^{n+1}, \surf^{n+1})} +  (\Delta t)^2  \, \int_{\Gamma_s^n} \omega\, (\bm u \cdot \bm n)^2\, ds + 2\,(\Delta t)^2 \, \int_{\Gamma_s^n} (\bm u \cdot \bm n)\, a^n\, ds \\
      &\leq  \|\surf^{n}\|_{L^2(\Omega_h^\perp)}^2 + 2\, \Delta t(a^n, \surf^n)_{\Omega^\perp_h} +(\Delta t)^2 \|a^n\|^2_{L^2(\Omega^\perp_h)} +  (\Delta t)^2\,\int_{\Gamma_s^n} \omega\, (\bm u \cdot \bm n)^2\, ds + 2\,(\Delta t)^2\, \int_{\Gamma_s^n} (\bm u \cdot \bm n)\, a^n\, ds.
  \end{split}
\end{equation}
We subtract $(\Delta t)^2  \, \int_{\Gamma_s^n} \omega\, (\bm u \cdot \bm n)^2\, ds + 2\,(\Delta t)^2 \, \int_{\Gamma_s^n} a^n\, (\bm u \cdot \bm n)\, ds$ on both sides of the inequality 
to finalize the proof of this proposition.
\end{proof}

\begin{proposition}[Stabilized explicit Euler domain volume conservation] \label{proposition:forward_euler_stabilized_domain_volume_cons}
  Solution to the explicit Euler scheme \eqref{eq:forward_euler_stabilized} with  source term $a^n=a(\bm x^\perp, t^n)$ 
  \reviewerTwo{and assumption $\surf^{n+1} > b$ everywhere on $\Omega_h^\perp$,} 
  conserves domain volume: $|\Omega_h^{n+1}| - |\Omega_h^n| = \Delta t\, \int_{\Omega^\perp_h} a^{n} \ud \bx^\perp$.
\end{proposition}
\begin{proof}
Applying all steps of the proof of Proposition \ref{proposition:backward_euler_domain_conservation} to \eqref{eq:forward_euler_stabilized}, proves this proposition.
\end{proof}

\subsection{Semi-implicit Euler discretized coupled problem \reviewerTwo{with the original FSSA stabilization and implicit surface regularization}}
\label{Sec:semi_fssa}
The semi-implicit Euler method stabilized using the FSSA approach \cite{Kaus_fssa, lofgren_fssa2} is formulated as follows. 
Given $\surf^n$ and $a^{n}$, find $\surf^{n+1} \in Z_h(\Omega_h^\perp)$, $\bm u^n \in \bm V_h(\Omega_h^n)$, and $\pi^n \in Q_h(\Omega_h^n)$, such that:
 \begin{equation}
  \label{eq:semi_backward_euler_fssa_stabilized}
  \left\{
\begin{aligned}
  \lft( \frac{\surf^{n+1} - \surf^n}{\Delta t}, w \rgt)_{\Omega_h^\perp} + (\bm u^{\perp,n} \cdot \nabla ^\perp \surf^{n+1} - u_z^{n}, w)_{\Omega_h^\perp} + \reviewerTwo{\edge(\surf^{n+1},w)}  &= (a^{n},w)_{\Omega_h^\perp}, \\
  2(\sqrt{\mu}\,\bm D\bm u^{n}, \sqrt{\mu}\,\bm D \bm v)_{\Omega_h^{n}} + \reviewerOne{\theta}\rho g\, \Delta t\, \int_{\Gamma_s^n} \!\!\!(\bm u \cdot \bm n)\, (\hat{\bz} \cdot \bm v)\, ds  - (\pi^{n}, \nabla \cdot  \bm v)_{\Omega_h^{n}} &= - \rho g (  \hat{\bz}, \bv)_{\Omega_h^{n}} 
  \reviewerTwo{- \reviewerOne{\theta} \rho g\, \Delta t\, \int_{\Gamma_s^n} \!\!\!(a \hat{\bz} \cdot \bm n)\, (\hat{\bz} \cdot \bm v)\, ds},\\
( \nabla \cdot  \bu^{n}, q)_{\Omega_h^{n}} &= 0,
\end{aligned}
\right.
\end{equation}
for all $\bm v \in \bm V_h(\Omega_h^n)$, $q \in Q_h(\Omega_h^n)$, and $w \in Z_h(\Omega_h^\perp)$. \reviewerOne{Here, $\theta$ is a user-defined variable, typically chosen between $0$ and $1$.}
\reviewerOne{The FSSA-stabilization was developed for geodynamics simulations \citep{Kaus_fssa} and recently adopted in} numerical models of \igorr{ice sheet dynamics \cite{lofgren_fssa1,lofgren_fssa2,tominec2024weakformshallowice}. In this setting the stabilization terms are typically scaled by a parameter $\theta$.} 
\igorr{The formulation can also be written in the explicit Euler setting, but numerical evidence suggests that 
the semi-implicit formulation allows taking \reviewerOne{slightly} larger stable time step sizes \cite{lofgren_fssa1,lofgren_fssa2}.}
\igorr{In both settings, the added non-symmetric term $\rho g\, \Delta t\, \int_{\Gamma_s^n} (\bm u \cdot \bm n)\, (\hat{\bz} \cdot \bm v)\, ds$ 
breaks the symmetry of the Stokes problem, which could be disadvantageous when using iterative linear solvers.} 

The stability analysis of \eqref{eq:semi_backward_euler_fssa_stabilized} requires extending the theoretical framework in the proof of Theorem \ref{theorem:surface_height_energy_rate} 
beyond the scope of the present paper. 
Instead, in Section \ref{sec:experiments:newtonian_tank} and Section \ref{sec:experiments:greenland} 
we experimentally show that the solution to this scheme: 
\begin{enumerate}
\item is conditionally stable with respect to $\Delta t$ in the total energy norm,
\item does not conserve the domain volume.
\end{enumerate}

\section{Numerical experiments}
\label{sec:experiments}
In the sections that follow, we outline numerical experiments verifying our theoretical results on the discrete stability of the Stokes free-surface flow. 
We use the open-source finite element library FEniCS, version 2019.1 \cite{Alnaes2014,Alnaes2015,UFL}.

\paragraph{Boundary conditions.} In all considered cases, the boundary conditions are:
$ \bm u = \bm 0 \text{ on } \Gamma_b$, $\bm u \cdot \bm n = 0 \text{ on } \Gamma_l$, and $\bm\sigma \cdot \bm n = 0 \text{ on } \Gamma_s.$
 
\paragraph{Mesh resolution.} We denote $\Delta x$ and $\Delta x^\perp$ as mesh sizes on $\Omega_h$ and $\Omega^\perp_h$ respectively. 
In some places, we report the results in terms of the number of mesh elements in horizontal direction $N_x$ and vertical direction $N_y$.

\paragraph{Stability criterion.} The stability criterion is based on \reviewerOne{the one‑step discrete analogue of Theorem \ref{theorem:surface_height_energy_rate}.}
\reviewerTwo{For the schemes considered below where the source term is evaluated explicitly at $t^n$ (EE-Unstabilized, EE‑Stabilized, EE‑FSSA, and SIE‑FSSA), 
we write the left and right-hand sides of the estimate as 
$E_L^n = \|\surf^{n+1}\|_{L^2(\Omega_h^\perp)}^2 + \frac{4\,\Delta t}{\rho g}\, \|\sqrt{\mu}\, \bm D \bm u\|^2_{L^2(\Omega_h^n)}$ and 
$E_R^n = \|\surf^{n}\|_{L^2(\Omega_h^\perp)}^2 + 2\, \Delta t(a^n, \surf^n)_{\Omega^\perp_h} + (\Delta t)^2 \|a^n\|^2_{L^2(\Omega^\perp_h)}$ respectively. Then 
we define the normalized energy at $t^n$ as:
\begin{equation}
  \label{eq:experiments_relative_normalized_energy}
\bar E^n = \frac{E_L^n - E_R^n}{\max_{n} |E_R^n|},\qquad n=0,1,2,...
\end{equation}
where $\bar E^n \leq 0$ implies stability per time step. }
\igorr{The above criterion is stated in the total energy norm. 
This norm reflects the coupling between $\bm{u}$ and $\surf$ and provides a natural measure of stability. 
In contrast, separate estimates for $\bm{u}$ and $\surf$ that ignore their interaction are not sufficient for analyzing the stability of the coupled system.}

\paragraph{Abbreviations.} 
\reviewerTwo{To distinguish between the different numerical schemes, we use the following abbreviations (all schemes include implicit surface regularization via $\edge(\surf^{n+1},w)$):} 
\emph{EE-Unstabilized} for \reviewerTwo{explicit Euler with implicit surface regularization (without coupled-energy stabilization)} \eqref{eq:forward_euler}, 
\emph{EE-Unstabilized-W} for EE-Unstabilized with $\Delta t$ computed using \eqref{eq:explicit_stability_final_weakerstab_tmp1} (weakly stable), 
\emph{EE-Stabilized} for \reviewerTwo{explicit Euler with coupled-energy stabilization and implicit surface regularization} \eqref{eq:forward_euler_stabilized}, 
\emph{SIE-FSSA} for the FSSA stabilized semi-implicit Euler scheme \eqref{eq:semi_backward_euler_fssa_stabilized}, 
\emph{EE-FSSA} for the FSSA stabilized explicit Euler scheme.

\subsection{Newtonian fluid dynamics in a tank}
\label{sec:experiments:newtonian_tank}
\begin{figure}[h!]
  \centering
\begin{tabular}{ccc}
  \multicolumn{3}{c}{\bf Newtonian fluid in a tank case}\vspace{0.1cm}  \\
  \textbf{Initial domain} & \textbf{Free-surface reference solution} & \\
  \includegraphics[width=0.34\linewidth]{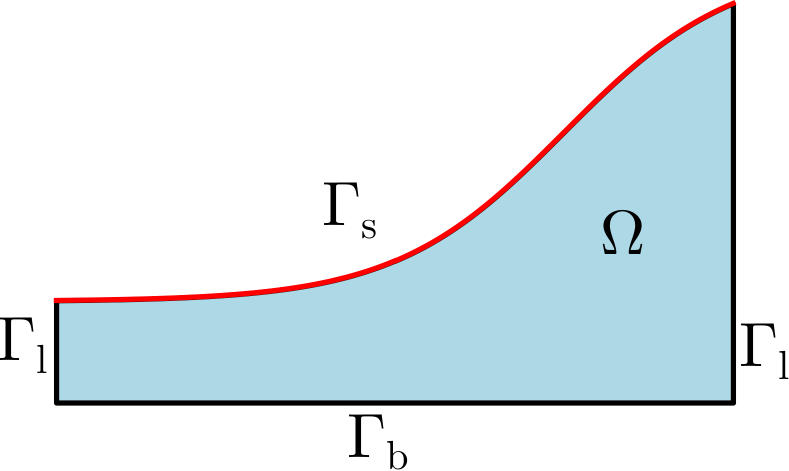} & 
  \raisebox{0.5cm}{\includegraphics[width=0.27\linewidth]{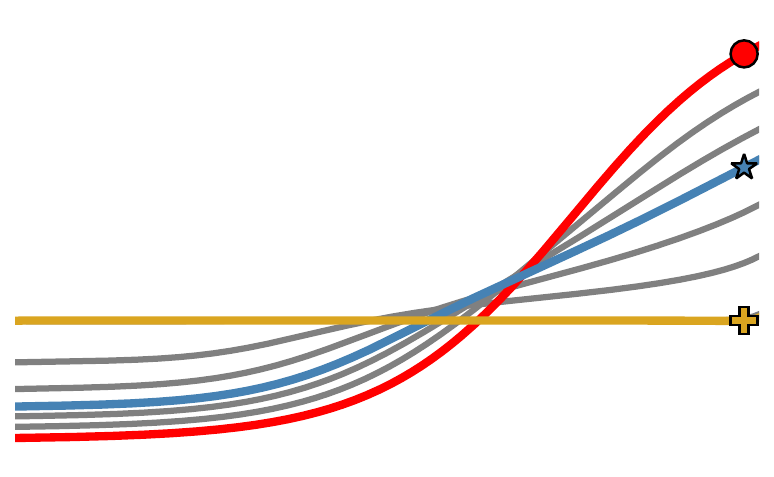}} &
  \raisebox{0.8cm}{\includegraphics[width=0.13\linewidth]{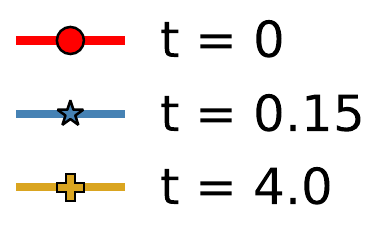}}
  \\
 (a) & (b) & \\
\end{tabular}
\caption{(a) Initial domain for Newtonian fluid in a tank case. (b) free-surface height evolving in time from $t=0$ to $t=4$.}
\label{fig:experiments:newtonian_tank:domain_solutions}
\end{figure}

In this section, we consider a Newtonian fluid case. This implies a linear Stokes problem. 
The domain $\Omega_h$ is defined as in \eqref{eq:surface_representations:Omega}, where we set $\Omega^\perp = [-1, 1]$, $b(x)=0$ and $\surf_{\reviewerTwo{0}}(x) = 0.5\, \tanh(2 x - 1) + 0.2$. 
In \eqref{eq:stokes_free_surface_weakform} we set $\mu = 0.3$ and $- \rho g \hat{\bz} = (0, -9.82)$. 
In Figure \ref{fig:experiments:newtonian_tank:domain_solutions} we display the domain shape and a reference solution \reviewerTwo{computed using 
\eqref{eq:forward_euler} under a very small time step and mesh size } with $a=0$ when the final simulation time is $\hat t = 4$. 

\subsubsection{Stability comparison of the time stabilization approaches.} 
We use $N_x = 120$ and $N_y = 120$ mesh elements in each direction and use the normalized energy \eqref{eq:experiments_relative_normalized_energy} as the stability criterion to verify the stability at different time steps when 
$\hat t = 4$. 
\begin{figure}[h!]
  \centering
\begin{tabular}{ccccc}
  \multicolumn{5}{c}{\bf Newtonian fluid in a tank} \\
  \multicolumn{5}{c}{\bf Normalized time-step energy difference} \vspace{0.1cm}\\
  \hspace{0.65cm}\textbf{EE-Unstabilized} & \hspace{0.5cm}\textbf{EE-Stabilized} & \hspace{0.4cm}\textbf{EE-FSSA} & \hspace{0.6cm}\textbf{SIE-FSSA}  \vspace{-0.1cm}\\
  \includegraphics[width=0.22\linewidth]{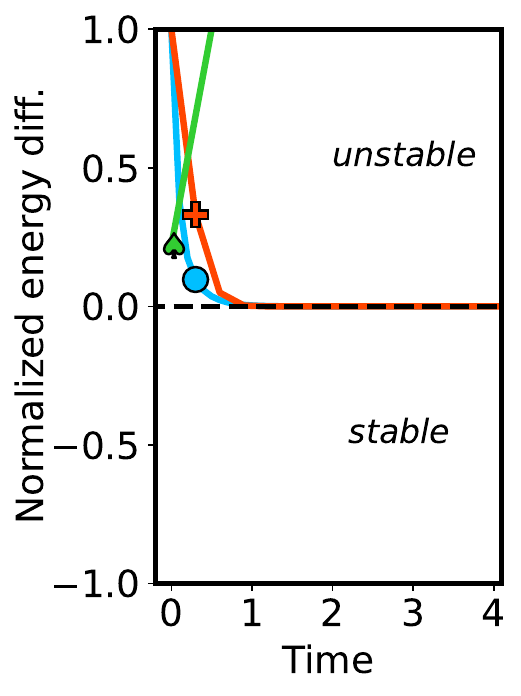} & 
  \hspace{-0.5cm}  \includegraphics[width=0.22\linewidth]{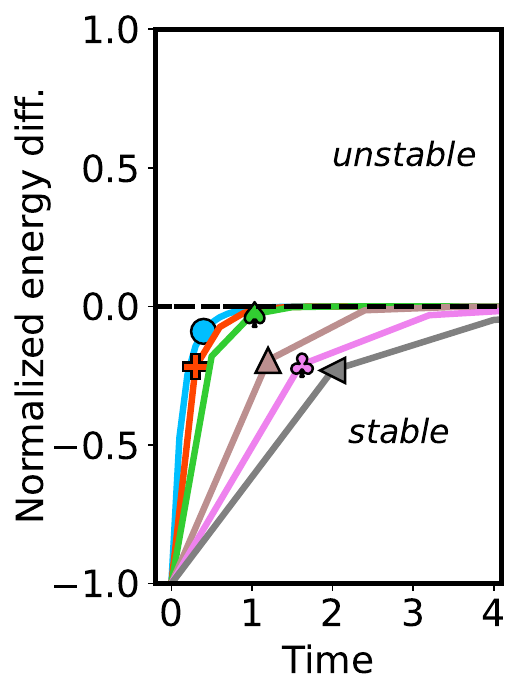} &  
  \hspace{-0.5cm}  \includegraphics[width=0.22\linewidth]{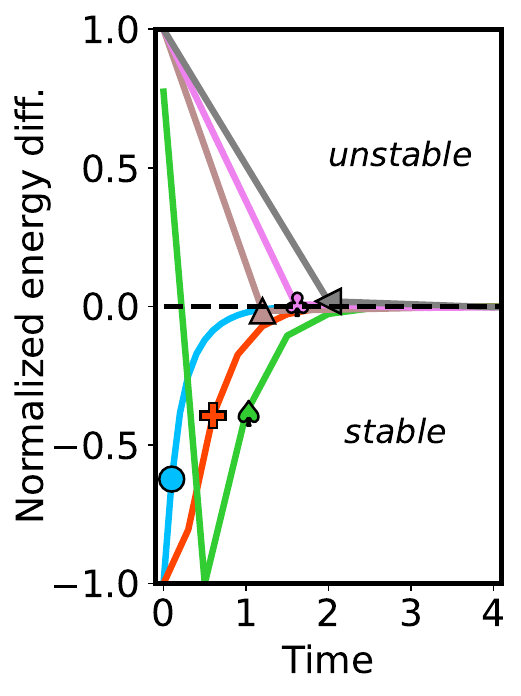} &  
  \hspace{-0.5cm}  \includegraphics[width=0.22\linewidth]{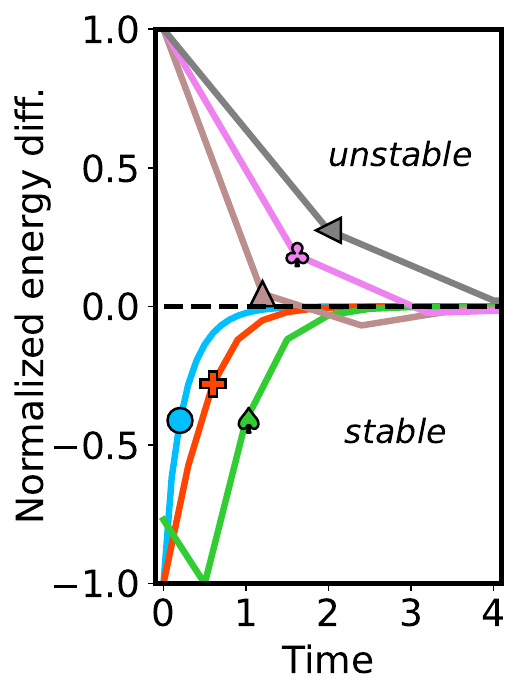} &  
  \hspace{-0.25cm}\raisebox{2cm}{\includegraphics[width=0.11\linewidth]{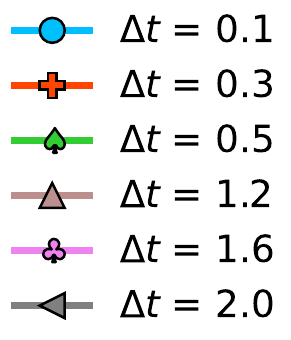}} \vspace{-0.1cm}\\
\end{tabular}
\caption{Normalized consecutive time step energy difference (stability criterion) for a Newtonian fluid case with source term $a=0$ 
compared across four different choices of time discretization.  Other parameters are: viscosity $\mu = 0.3$, uniform mesh $N_x=N_y=120$, final time $\hat t = 4$.}
\label{fig:experiments:newtonian_tank:consecutive_energy_comparison}
\end{figure}
Results are given in Figure \ref{fig:experiments:newtonian_tank:consecutive_energy_comparison}. 
First, we observe that the EE-Unstabilized scheme is unstable for all the considered time steps, confirming our theoretical observation provided in 
Proposition \ref{proposition:stability_forward_euler}. 
The edge stabilization term \reviewerTwo{$\edge(\surf^{n+1},w)$} \eqref{eq:discrete_problem_analysis:edge_stabilization} on itself does not correct the stability in time. 
Second, we observe that the EE-Stabilized scheme, with our new stabilization approach, is always stable, which confirms our theoretical 
result in Proposition \ref{proposition:stability_forward_euler_stabilized}. 
Third, the EE-FSSA scheme and the SIE-FSSA scheme are only stable for some time steps. 

\begin{figure}[h!]
  \centering
\begin{tabular}{ccccc}
  \multicolumn{5}{c}{\bf Newtonian fluid in a tank with source term} \\
  \multicolumn{5}{c}{\bf Normalized time-step energy difference} \vspace{0.1cm}\\
  \hspace{0.65cm}\textbf{EE-Unstabilized} & \hspace{0.5cm}\textbf{EE-Stabilized} & \hspace{0.5cm}\textbf{EE-FSSA} & \hspace{0.6cm}\textbf{SIE-FSSA}  \vspace{-0.1cm}\\
  \includegraphics[width=0.22\linewidth]{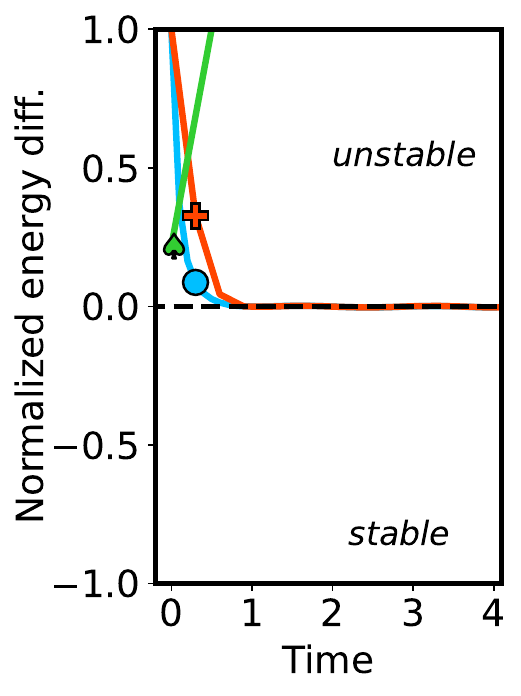} & 
  \hspace{-0.5cm}  \includegraphics[width=0.22\linewidth]{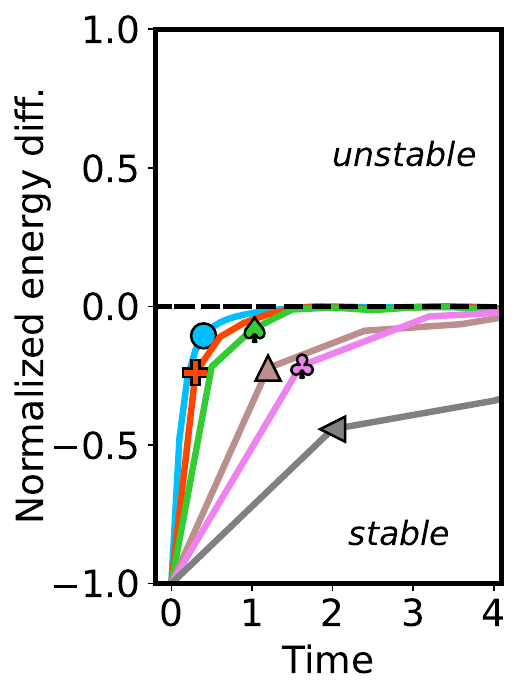} &  
  \hspace{-0.5cm}  \includegraphics[width=0.22\linewidth]{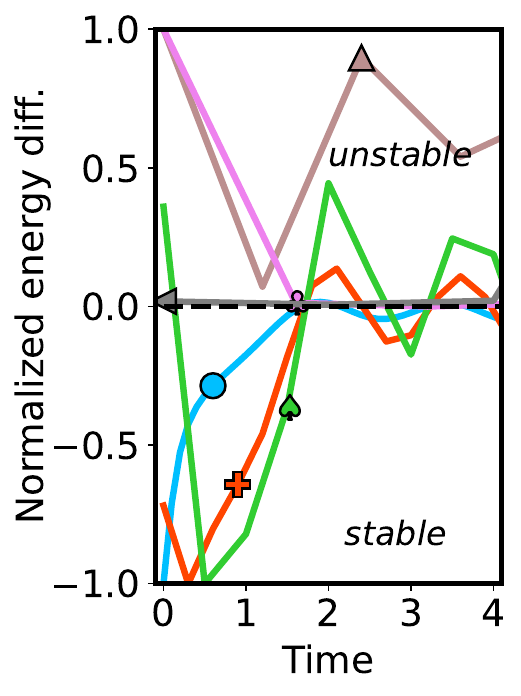} &  
  \hspace{-0.5cm}  \includegraphics[width=0.22\linewidth]{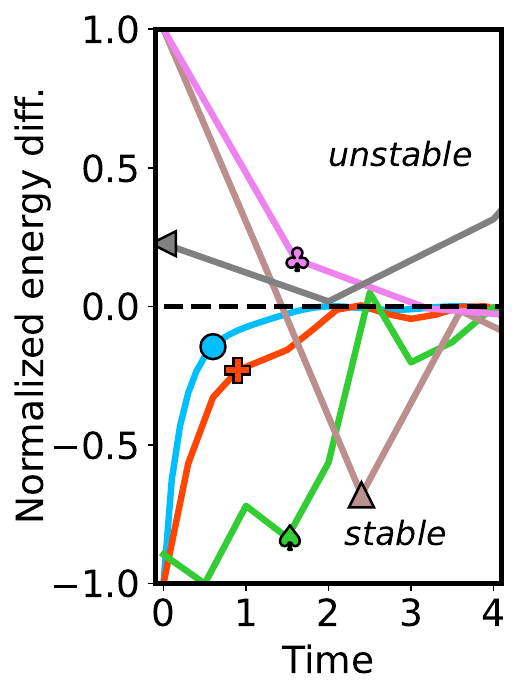} &  
  \hspace{-0.25cm}\raisebox{2cm}{\includegraphics[width=0.11\linewidth]{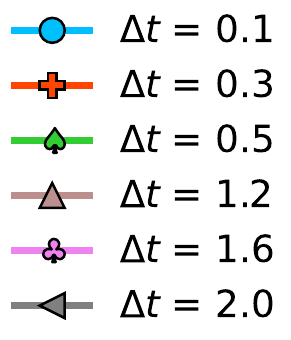}} \vspace{-0.1cm}\\
\end{tabular}
\caption{Normalized consecutive time step energy difference (stability criterion) for a Newtonian fluid case with source term $a= 0.2\, (x^\perp)^2\, (0.3 + \sin(x^\perp))\sin(2t)\, \, $ 
compared across four different choices of time discretization.  Other parameters are: viscosity $\mu = 0.3$, uniform mesh $N_x=N_y=120$, final time $\hat t = 4$.}
\label{fig:experiments:newtonian_tank:consecutive_energy_comparison_mass_balance}
\end{figure}

We also test the numerical stability when the source term takes values $a=a(\bm x^\perp, t) = 0.2\, (x^\perp)^2\, (0.3 + \sin(x^\perp))\sin(2t)\, \, $.
The result is given in Figure \ref{fig:experiments:newtonian_tank:consecutive_energy_comparison_mass_balance}. 
We find the interpretation of these results analogous to those when $a=0$. 
\reviewerTwo{Hence, our stability estimates under a non-zero source term are confirmed. Moreover, 
the source term tailored stabilization in EE-Stabilized is numerically validated.}

Overall, we find the newly proposed EE-Stabilized scheme \eqref{eq:forward_euler_stabilized} best in terms of stability as it allows arbitrarily large time step sizes. 

\subsubsection{Numerical errors under mesh and time step size refinement for different stabilization approaches.} 
\igorr{We numerically test the convergence of the numerical errors in $\bm u$, $\bm u^\perp$, and $\surf$ towards $0$. 
This is only possible when simultaneously refining $\Delta x \to 0$, $\Delta x^\perp \to 0$, and $\Delta t\, \to 0$}. 
We examine how adding different stabilization terms affects the numerical error, i.e., to what extent the stabilizations are consistent with the discretized problem. 
In all the tested cases we simultaneously refine the mesh size and the time step size according to Table \ref{table:experiments:newtonian:convergence_error:meshsize_timestep}. 
We chose the time step sizes so that they satisfy the time step restriction \eqref{eq:explicit_stability_final_weakerstab_tmp1}, 
required when using EE-Unstabilized-W.
\begin{figure}
    \centering
 {\bf Newtonian fluid in a tank\\ Table of time steps and the corresponding mesh sizes for the convergence test} \vspace{0.1cm}\\
  \begin{tabular}{|c|c|c|c|c|}
    \hline
  $\Delta t$ & $0.5$ & $0.25$ & $0.125$ & $0.0625$ \\
  \hline
  $\Delta x^\perp$ & $0.2$ & $0.1$ & $0.05$ & $0.025$\\
  \hline
  $\Delta x$ & $0.24$ & $0.11$ & $0.06$ & $0.03$ \\
  \hline
  \end{tabular}\vspace{0.2cm}
  \caption{Time step $\Delta t$ and the corresponding mesh sizes $\Delta x^\perp$ on $\Omega_h^\perp$ and $\Delta x$ on $\Omega_h$ 
 for the error convergence test in the Newtonian fluid case.}
  \label{table:experiments:newtonian:convergence_error:meshsize_timestep}
\end{figure}
The final simulation time is $\hat t = 1$. We computed the relative numerical errors against a numerical reference 
using the EE-Unstabilized-W scheme with mesh sizes $\Delta x = 0.015$, $\Delta x^\perp = 0.0125$, and time step $\Delta t = 0.005$. 
\begin{figure}[h!]
  \centering
\begin{tabular}{cccc}
  \multicolumn{4}{c}{\bf Newtonian fluid in a tank} \\
  \multicolumn{4}{c}{\bf Convergence of error under simultaneous mesh and time refinement} \vspace{0.1cm}\\
  \hspace{1.2cm}\textbf{Surface height} & \hspace{0.5cm}\textbf{Surface velocity} & \hspace{0.6cm}\textbf{Velocity}  \vspace{-0.05cm}\\
  \hspace{0.1cm}\includegraphics[width=0.23\linewidth]{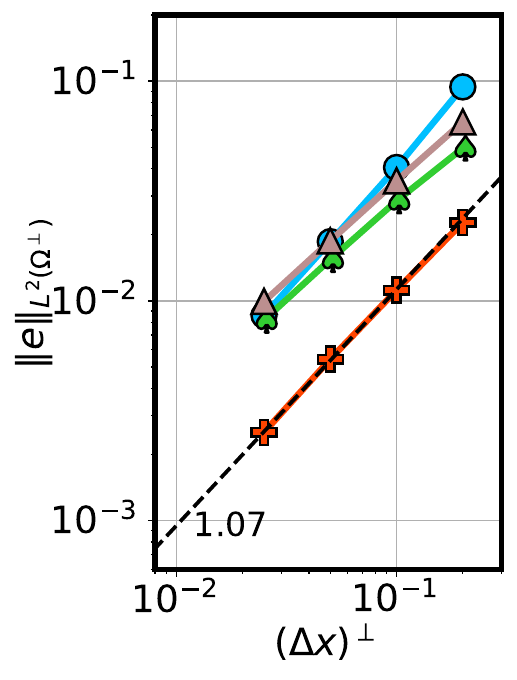} & 
  \hspace{-0.3cm}\includegraphics[width=0.23\linewidth]{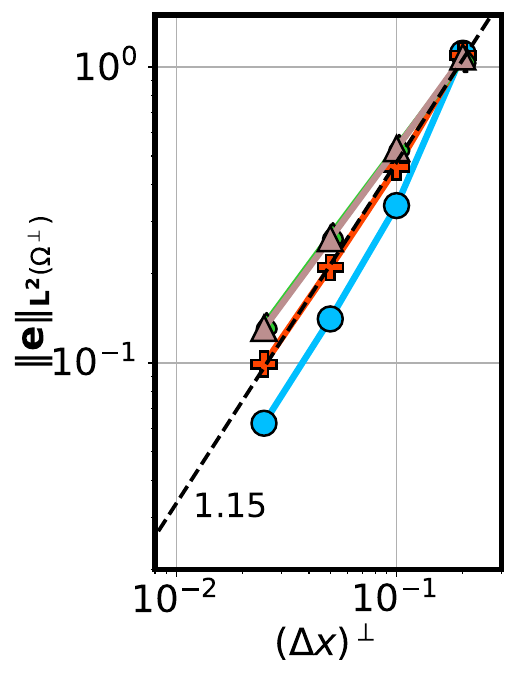} & 
  \hspace{-0.3cm}\includegraphics[width=0.23\linewidth]{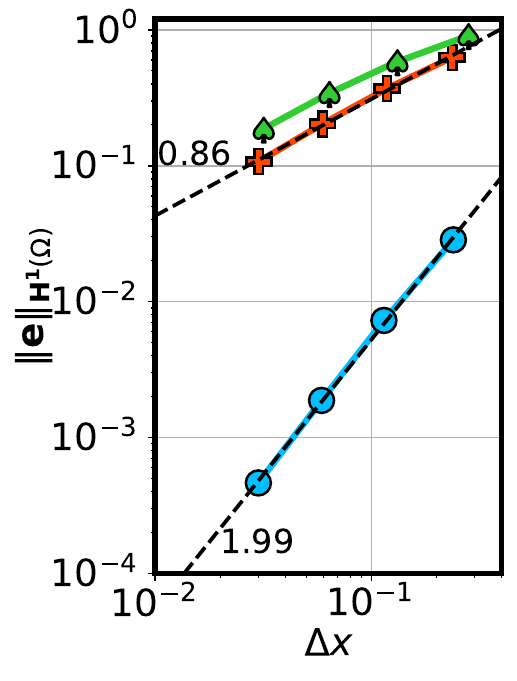} & 
  \hspace{-0.25cm}\raisebox{2.7cm}{\includegraphics[width=0.2\linewidth]{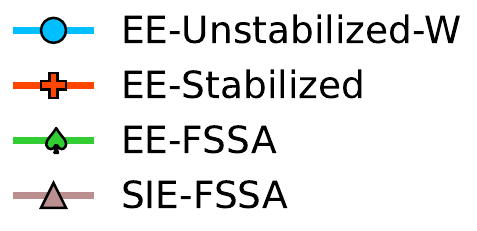}} \vspace{-0.1cm}\\
 
  \hspace{1.0cm}\includegraphics[width=0.155\linewidth]{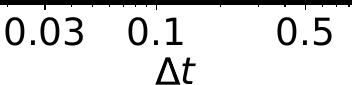} & 
  \hspace{0.55cm}\includegraphics[width=0.155\linewidth]{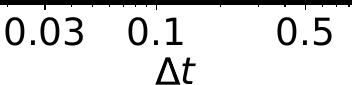} & 
  \hspace{0.57cm}\includegraphics[width=0.155\linewidth]{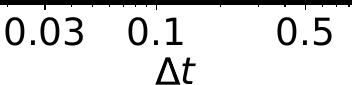} &

\end{tabular}
\caption{Convergence of the relative numerical error (Newtonian fluid case) under simultaneous mesh and time refinement in surface height on $\Omega^\perp$, surface velocity on $\Omega^\perp$, and velocity on $\Omega$, 
for different choices of discretization in time and the corresponding stabilization approaches. 
Parameters in all cases are: viscosity coefficient $\mu=0.3$, source term $a=0$.}
\label{fig:experiments:newtonian_tank:error_convergence}
\end{figure}

In Figure \ref{fig:experiments:newtonian_tank:error_convergence} we 
report the numerical errors in the surface height function on $\Omega^\perp$, surface velocities on $\Omega^\perp$, and the velocities on $\Omega$. 
\reviewerTwo{Our intention is both to understand whether the stabilization terms are consistent with the discretized problem and to compare the errors across the different schemes.}
We observe that the error in the surface height case and the surface velocity case converges with order $1$ for all the considered schemes, which is expected. 
The EE-Stabilized scheme is significantly more accurate compared to the EE-Unstabilized-W scheme. 

In the surface velocity case, the error is comparable across all the schemes. 

In the velocity case, the error decays with order $2$ when the EE-Unstabilized-W scheme is used, but the error decay is of order $1$ in all other cases. This implies 
that the stabilization term added to the Stokes problem is consistent up to the first order in each of the cases EE-Stabilized, SIE-FSSA, 
and EE-FSSA. The velocity error is significantly 
smaller in the EE-Unstabilized-W case compared to all the stabilized cases. Among those, the error is smallest in the EE-Stabilized case \reviewerTwo{in $\surf$}.

Among all the stabilized schemes, the EE-Stabilized scheme induced the smallest errors. The main point of interest in the free-surface flow simulations is the surface height, where 
EE-Stabilized was, in addition, more accurate than EE-Unstabilized-W. 
\reviewerTwo{Overall, the stabilization terms added in the EE-Stabilized scheme are consistent with the discretized problem and did not introduce significant errors 
as long as the main point of interest is the surface height.}

\subsubsection{Conservation of domain volume}
\begin{figure}[h!]
  \centering
\begin{tabular}{cccc}
  \multicolumn{4}{c}{\bf Newtonian fluid in a tank: Domain volume error in time} \vspace{0.05cm}\\
  \hspace{1cm}\textbf{EE-Stabilized} & \hspace{0.6cm}\textbf{EE-FSSA} & \hspace{0.6cm}\textbf{SIE-FSSA}   \vspace{-0.1cm}\\
  \includegraphics[width=0.25\linewidth]{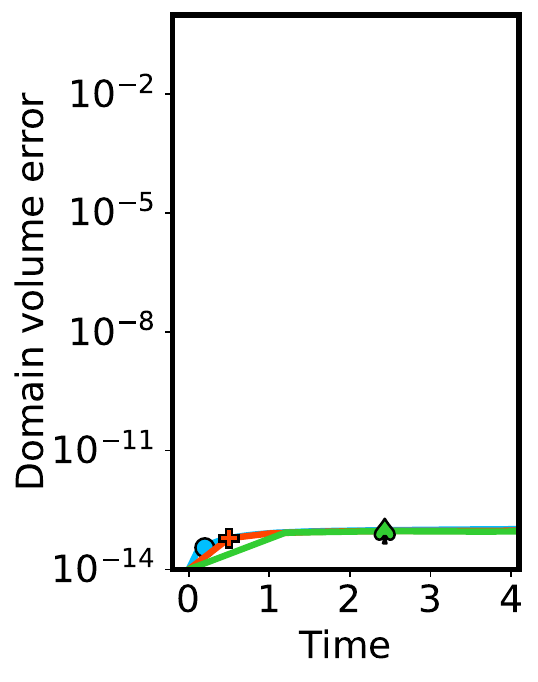} & 
  \hspace{-0.5cm}  \includegraphics[width=0.25\linewidth]{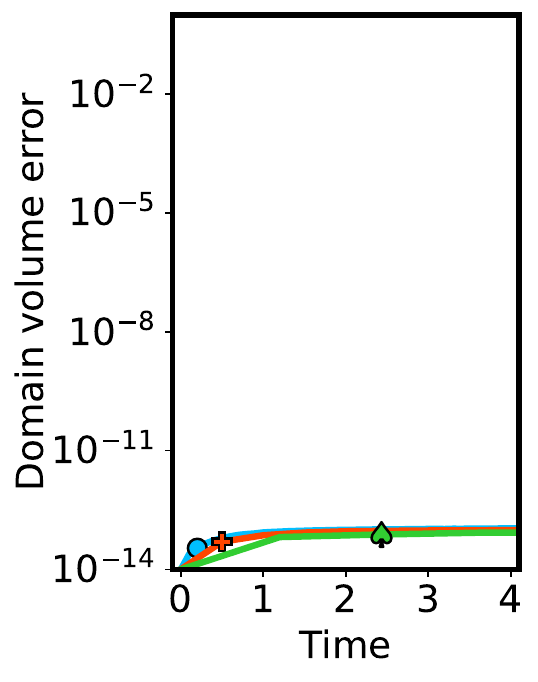} &  
  \hspace{-0.5cm}  \includegraphics[width=0.25\linewidth]{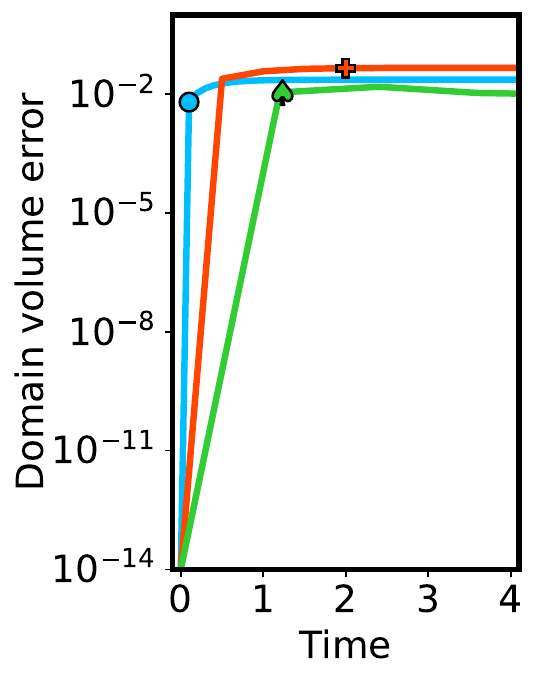} &  
  \hspace{-0.25cm}\raisebox{3.1cm}{\includegraphics[width=0.13\linewidth]{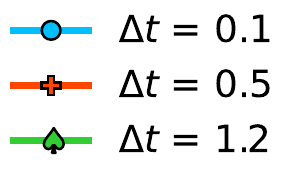}} \vspace{-0.1cm}\\
\end{tabular}
\caption{Domain $\Omega$ volume error (Newtonian fluid case) as a function of time for different choices of discretization in time and the corresponding stabilizations. 
We use a uniform mesh $N_x = N_y = 120$, viscosity $\mu = 0.3$, source term $a=0$.}
\label{fig:experiments:newtonian_tank:mass_in_time}
\end{figure}

Ideally, a numerical scheme for solving the free-surface coupled Stokes problem conserves the domain volume in time, i.e., the domain volume error is at the level of round-off errors (close to machine precision). 
In this experiment, we use $N_x = 120$ and $N_y = 120$ mesh elements in each direction and set the final simulation time to $\hat t = 4$. 
In Figure \ref{fig:experiments:newtonian_tank:mass_in_time}, we display the domain volume error for different time step sizes. In the EE-Stabilized case, we observe that the domain volume is conserved up to round-off errors, 
which is in line with Proposition \ref{proposition:forward_euler_domain_volume_cons}. 
The EE-FSSA scheme also conserves the domain volume. 
On the other hand, the SIE-FSSA scheme does not conserve the domain volume as observed from the figure.
\subsection{Motivation for suppressing high-frequency oscillations due to unresolved gradients.} 
\begin{figure}[h!]
  \centering
\begin{tabular}{ccc}
  \multicolumn{2}{c}{\bf Newtonian fluid in a tank: surface height oscillations} \vspace{0.05cm}\\
  \hspace{0.65cm}\textbf{Unstabilized} & \hspace{0.5cm}\textbf{Edge stabilized} & \vspace{-0.05cm}\\
  \includegraphics[width=0.34\linewidth]{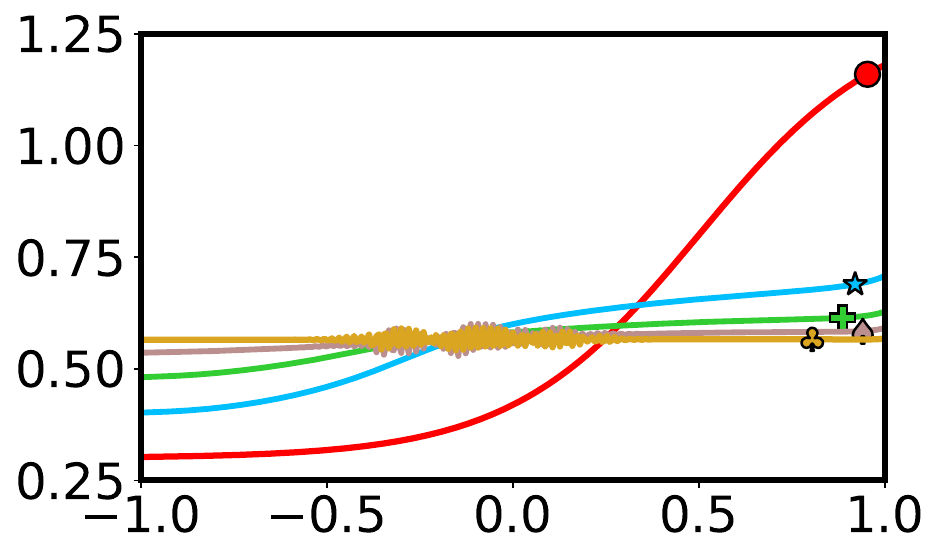} &
  \includegraphics[width=0.34\linewidth]{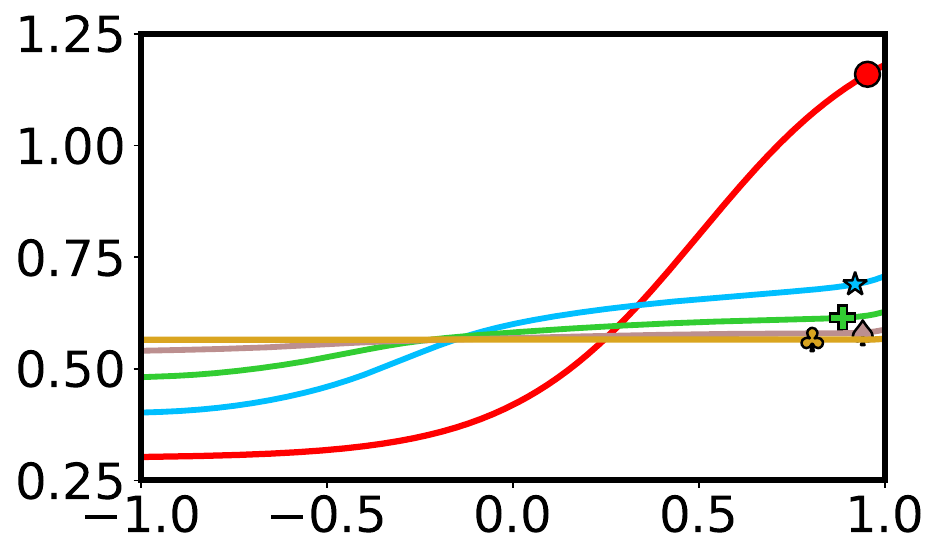} &
  \raisebox{0.8cm}{\includegraphics[width=0.12\linewidth]{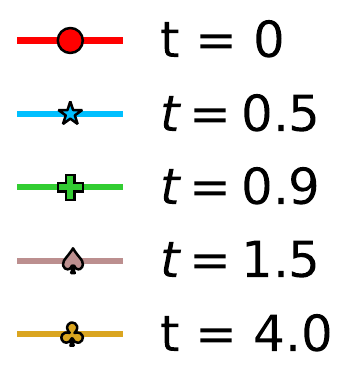}}
\end{tabular}
\caption{Surface height solutions (Newtonian fluid) at different $t$ in the EE-Stabilized case without additional stabilization of the high-frequency oscillations (left plot), 
and with stabilization of the high-frequency oscillations using the edge stabilization term \reviewerTwo{$\edge(\surf^{n+1},w)$} \eqref{eq:discrete_problem_analysis:edge_stabilization}. 
We use a uniform mesh $N_x = N_y = 300$, viscosity $\mu = 0.3$, source term $a=0$.}
\label{fig:experiments:newtonian_tank:gibbs}
\end{figure}

Even when the solution is stable in accordance with the stability estimate stated in 
\reviewerTwo{Proposition} \ref{proposition:stability_forward_euler}, the free-surface function can in practice entail high-frequency oscillations, particularly when the mesh resolution is large. Theoretical reasoning is provided in the scope of 
Section \ref{sec:continuous_problem_analysis:well_posedness_discussion} and Section \ref{sec:discrete_problem_analysis:continuity_regularization}. 
To suppress the oscillations, we use the edge stabilization term \reviewerTwo{$\edge(\surf^{n+1},w)$} \eqref{eq:discrete_problem_analysis:edge_stabilization}. 
In Figure \ref{fig:experiments:newtonian_tank:gibbs} we visualize the free-surface solution computed from the \reviewerTwo{explicit Euler scheme with coupled-energy stabilization} \eqref{eq:forward_euler_stabilized} 
at different times with a small time step ($\Delta t = 0.05$) and a fine mesh ($N_x=N_y=300$), with and without
edge stabilization. We observe that the edge stabilization term \reviewerTwo{$\edge(\surf^{n+1},w)$} effectively filters out high-frequency oscillations.

\subsection{\reviewerTwo{Non-Newtonian fluid dynamics over \reviewerTwo{a thin and rough geometry}}}
\label{sec:experiments:greenland}
In this section, we experimentally verify that our theory is correct for non-Newtonian fluids (nonlinear Stokes problem) on 
\reviewerTwo{a thin domain} with a rough boundary shape.

\begin{figure}[h!]
  \centering
\begin{tabular}{ccc}
  \multicolumn{3}{c}{\bf \reviewerTwo{Non-Newtonian fluid dynamics over a thin and rough geometry}}\vspace{0.1cm}  \\
  \textbf{Initial domain} & \textbf{Free-surface reference solution} & \vspace{0.05cm}\\
  \includegraphics[width=0.34\linewidth]{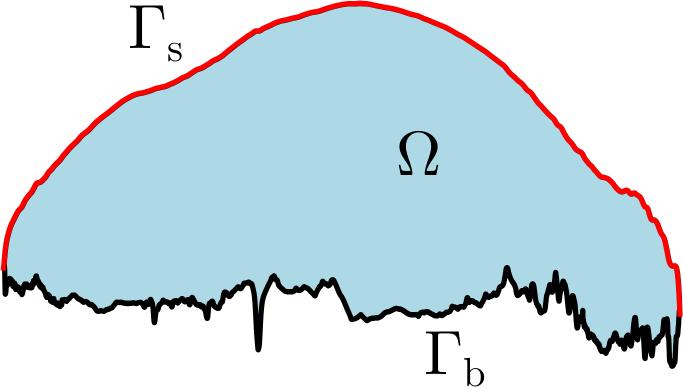} & 
  \raisebox{0.4cm}{\includegraphics[width=0.34\linewidth]{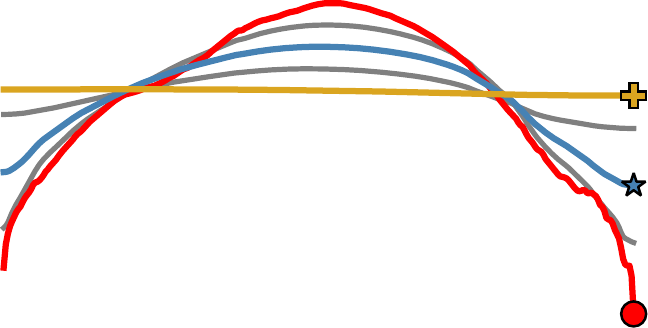}} &
  \raisebox{0.8cm}{\includegraphics[width=0.15\linewidth]{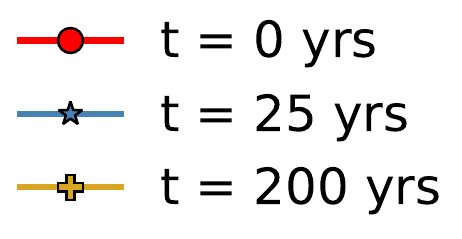}}
  \\
 (a) & (b) & \\
\end{tabular}
\caption{\reviewerTwo{(a) Initial domain for the non-Newtonian fluid case. (b) free-surface height evolving in time from $t=0$ years to $t=200$ years.}}
\label{fig:experiments:greenland:domain_solutions}
\end{figure}

The domain $\Omega_h$ as illustrated in Figure \ref{fig:icesheet_sketch} was extracted from the 3D Greenland data set BedMachine \cite{bedmachine} and 
interpolated using cubic splines to form the functions $\surf_0(\bm x^\perp)$ and $b(\bm x^\perp)$.  \reviewerTwo{The simulation is, however, not intended as a realistic ice sheet simulation and the evolution of the surface is not a representation of ice sheet dynamics.  
We omit lateral traction, set $a=0$, and use a short $200$-year window solely to stress-test numerical stability and volume conservation.}
In our simulations, we use ice parameters in SI units. 
We set $\Omega^\perp = [-428675 \text{ m}, 489475 \text{ m}]$. 
In \eqref{eq:stokes_free_surface_weakform} we set density $\rho =  910 \, \text{kg} \, \text{m}^{-3}$, gravity $g = 9.82 
\, \text{m} \, \text{s}^{-2}$. 
In \eqref{eq:stokes_viscosity} we set $p=4/3$ and viscosity coefficient (ice softness parameter) 
$\mu_0=3.1688 \cdot 10^{-24} \, (\text{Pa})^{-3} \text{s}^{-1}$. 
In all cases, we use $N_x = 300$ and $N_y = 20$ mesh elements in horizontal and vertical directions, consider a homogeneous source term $a=0$, 
use the Picard iteration to solve the nonlinear Stokes problem with a relative convergence tolerance $10^{-6}$, 
and perform the simulation until final time $\hat t = 200 \text{ years}$. 
In Figure \ref{fig:experiments:greenland:domain_solutions} we display the domain shape and a reference solution. 
\reviewerTwo{The reference solution is computed using the EE-Stabilized scheme with a very small time step and mesh size. 
The present experiment is intended to verify our theory in the non-Newtonian fluid case over a rough and thin geometry.}

\subsubsection{Stability comparison of the time stabilization approaches.} 
\begin{figure}[h!]
  \centering
\begin{tabular}{ccccc}
  \multicolumn{5}{c}{\bf \reviewerTwo{Non-Newtonian fluid dynamics over a thin and rough geometry}} \\
  \multicolumn{5}{c}{\bf Normalized time-step energy difference} \vspace{0.1cm}\\
  \hspace{0.65cm}\textbf{EE-Unstabilized} & \hspace{0.3cm}\textbf{EE-Stabilized} & \hspace{0.3cm}\textbf{EE-FSSA} & \hspace{0.4cm}\textbf{SIE-FSSA}  \vspace{-0.1cm}\\
  \includegraphics[width=0.21\linewidth]{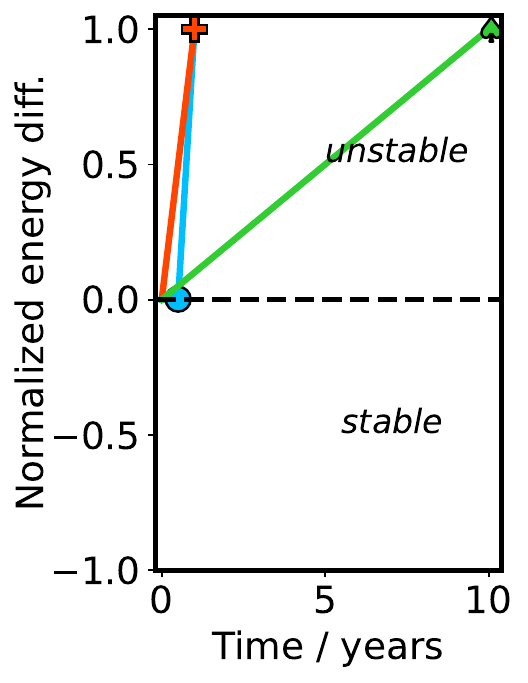} & 
  \hspace{-0.5cm}  \includegraphics[width=0.22\linewidth]{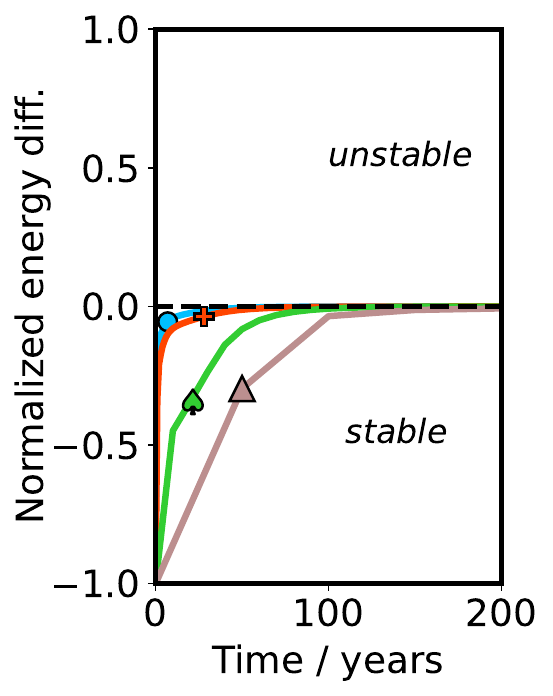} &  
  \hspace{-0.5cm}  \includegraphics[width=0.22\linewidth]{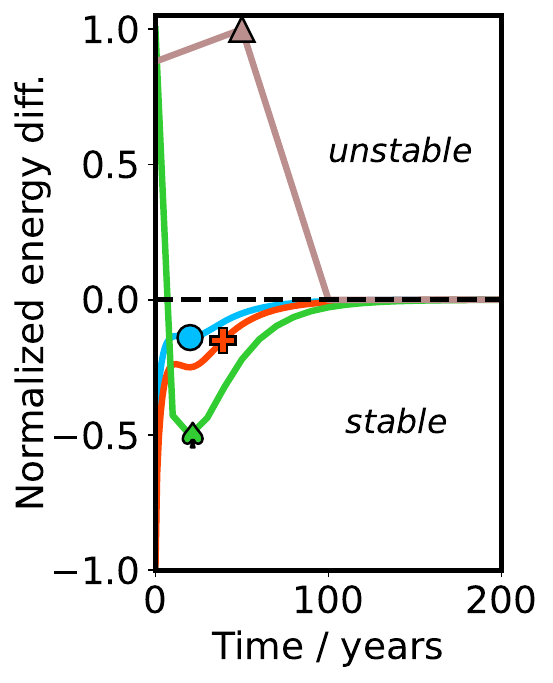} &  
  \hspace{-0.5cm}  \includegraphics[width=0.22\linewidth]{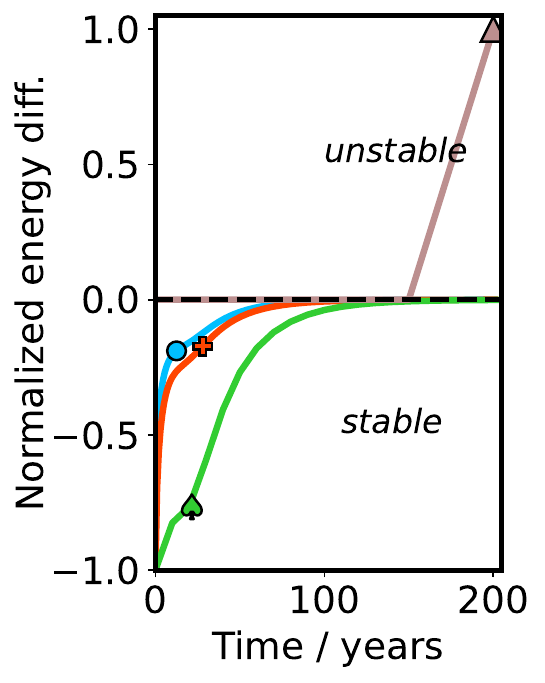} &  
  \hspace{-0.5cm}\raisebox{2.5cm}{\includegraphics[width=0.14\linewidth]{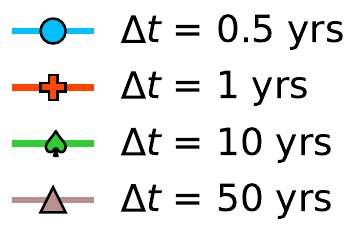}} \vspace{-0.1cm}\\
\end{tabular}
\caption{Normalized consecutive time step energy difference (stability criterion) for non-Newtonian fluid case with source term $a=0$ 
compared across four different choices of time discretization.  We use a non-uniform mesh with $N_x= 300$ elements in the horizontal direction 
and $N_y = 20$ elements in the vertical direction.}
\label{fig:experiments:greenland:consecutive_energy_comparison}
\end{figure}

To examine numerical stability we measure the normalized energy difference \eqref{eq:experiments_relative_normalized_energy} 
in each time step. Results are given in Figure \ref{fig:experiments:greenland:consecutive_energy_comparison}. 
First, the EE-Unstabilized scheme is unstable for all the considered time steps, confirming 
the result in Proposition \ref{proposition:stability_forward_euler}. 
Second, EE-FSSA and SIE-FSSA schemes are unstable when using the largest considered time step $\Delta t = 50$ years. 
The edge stabilization term \reviewerTwo{$\edge(\surf^{n+1},w)$} \eqref{eq:discrete_problem_analysis:edge_stabilization} does not prevent instabilities in time. 
Third, the EE-Stabilized scheme is stable for all the considered time steps, which confirms our theoretical result in 
Proposition \ref{proposition:stability_forward_euler_stabilized}.

\begin{figure}[h!]
  \centering
  \begin{tabular}{ccc}
    \multicolumn{3}{c}{\bf \reviewerTwo{Non-Newtonian fluid dynamics over a thin and rough geometry}} \\
    \multicolumn{3}{c}{\bf Free-surface height solutions with different stabilizations in time} \vspace{0.1cm}\\
    \hspace{0.5cm}\textbf{EE-Unstabilized} & \textbf{EE-Stabilized} & \textbf{SIE-FSSA}  \vspace{0.1cm}\\
    \includegraphics[width=0.3\linewidth]{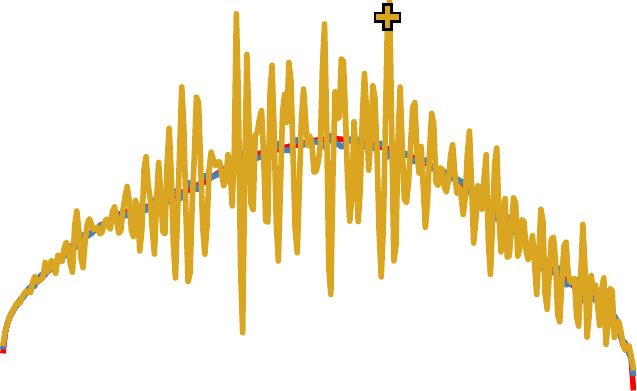} & 
    \hspace{-0.1cm} \includegraphics[width=0.3\linewidth]{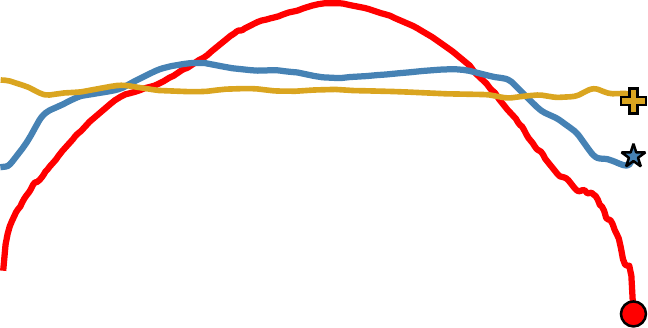} &  
    \hspace{-0.1cm}   \includegraphics[width=0.3\linewidth]{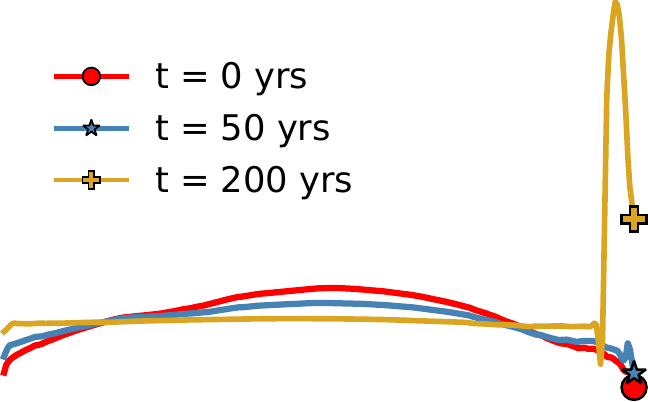}
  \end{tabular}
  \caption{Numerical free-surface height solutions varying in time, when using schemes: EE-Unstabilized with $\Delta t = 0.5$ years, EE-Stabilized with $\Delta t = 50$ years, 
 and SIE-FSSA with $\Delta t = 50$ years. 
 Edge stabilization term \reviewerTwo{$\edge(\surf^{n+1},w)$} is added to all of the schemes. We use $N_x=300$ mesh elements in the horizontal direction and $N_y=20$ mesh elements in the vertical direction.}
  \label{fig:experiments:greenland:numerical_solutions_diff_stabs}
\end{figure}

To visually inspect the stability behaviour of the numerical solutions, we plot the free-surface height solutions in Figure \ref{fig:experiments:greenland:numerical_solutions_diff_stabs}, 
for all the considered numerical schemes. 
In the EE-Unstabilized and SIE-FSSA cases, the chosen time step corresponds to the minimal value under which the solution becomes unstable. 
In the figure we observe that although the EE-Stabilized solution is inaccurate at $t=200$ years due to a large time step $\Delta t = 50$ years 
when compared to the reference solution in Figure \ref{fig:experiments:greenland:domain_solutions}, 
it is the only example of a solution that does not entail instabilities.

\subsubsection{Conservation of domain volume}
\begin{figure}[h!]
  \centering
\begin{tabular}{ccccc}
  \multicolumn{4}{c}{\bf \reviewerTwo{Non-Newtonian fluid dynamics over a thin and rough geometry}} \\
  \multicolumn{4}{c}{\bf Relative domain volume error in time} \vspace{0.1cm} \\
  \hspace{0.8cm}\textbf{EE-Stabilized} & \hspace{0.4cm}\textbf{EE-FSSA} & \hspace{0.4cm}\textbf{SIE-FSSA} & \vspace{-0.13cm}\\
  \includegraphics[width=0.25\linewidth]{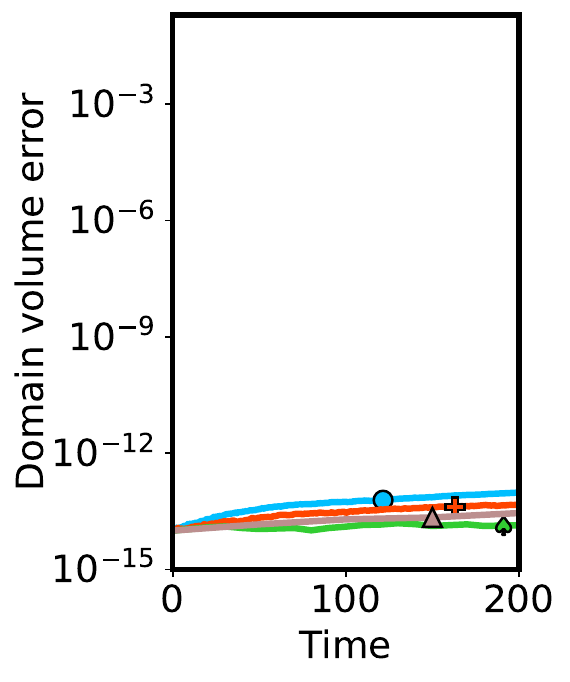} & 
  \hspace{-0.5cm}  \includegraphics[width=0.25\linewidth]{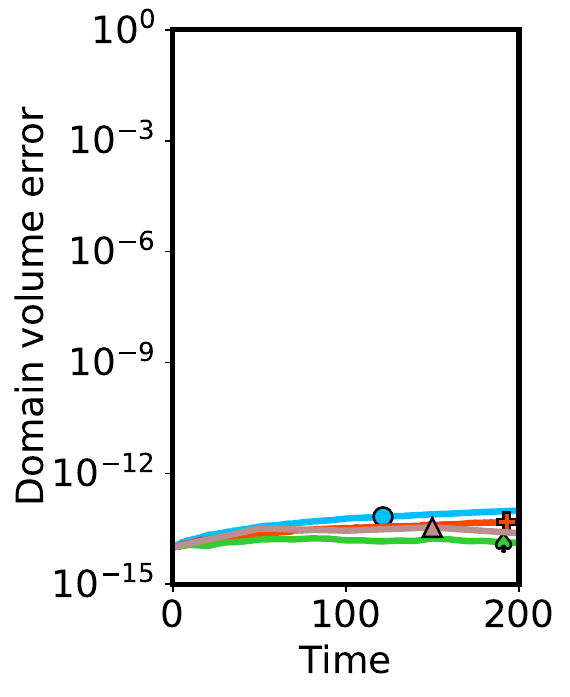} &  
  \hspace{-0.5cm}  \includegraphics[width=0.25\linewidth]{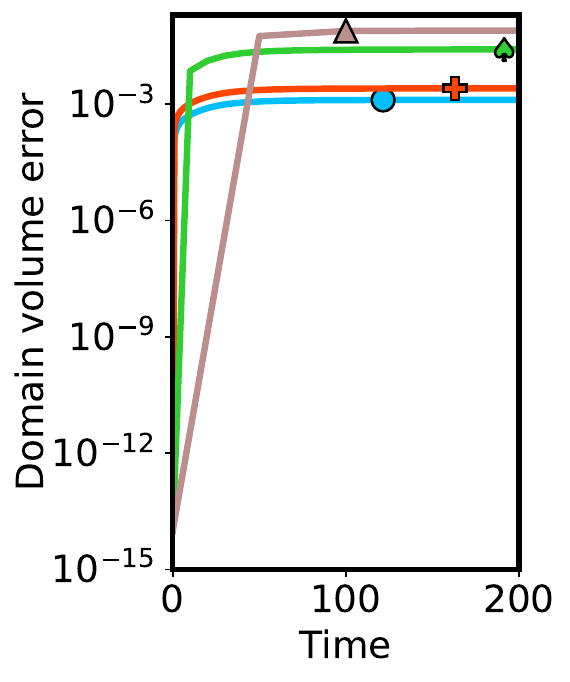} &  
  \hspace{-0.25cm}\raisebox{2.7cm}{\includegraphics[width=0.14\linewidth]{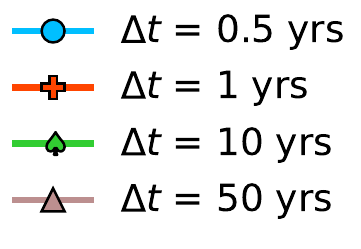}} \vspace{-0.1cm}\\
\end{tabular}
\caption{Domain $\Omega$ volume error (non-Newtonian fluid case) as a function of time for different choices of discretization in time and the corresponding stabilizations. 
We use a non-uniform mesh with $N_x = 300$ elements in the horizontal direction and $N_y = 20$ elements in the vertical direction.}
\label{fig:experiments:greenland:mass_in_time}
\end{figure}

In Figure \ref{fig:experiments:greenland:mass_in_time}, we display the relative domain volume error for different time step sizes. In the EE-Stabilized case, we observe that the domain volume is conserved up to round-off errors, 
which confirms Proposition \ref{proposition:forward_euler_domain_volume_cons}. 
The EE-FSSA scheme also conserves the domain volume. 
On the other hand, the SIE-FSSA scheme does not conserve the domain volume as observed from the figure.

\section{Final remarks}
\label{sec:final_remarks}

We have investigated the stability and conservation properties of the coupled Stokes/free-surface system. 
The analysis covered Newtonian (linear Stokes) and non-Newtonian (nonlinear Stokes) fluids, both at the continuous and fully discrete levels. 
\reviewerTwo{Our analysis gives $L^2$-norm stability estimates which are 
an important building block when proving well-posedness. 
These estimates are a necessary building block toward well-posedness once additional regularization ensures Lipschitz continuity 
of the surface height function, see Section \ref{sec:continuous_problem_analysis:well_posedness_discussion}. 

To practically remedy the lack of default regularity of $\surf$, 
we employed an edge stabilization term \reviewerTwo{$\edge(\surf^{n+1},w)$} \eqref{eq:discrete_problem_analysis:edge_stabilization}. 
In all of the numerical schemes the term is treated implicitly so as to not change the time-step size restrictions. 
The implicit treatment does not significantly increase the computational cost of the Lagrange FEM-based explicit Euler scheme that 
in any case requires solving a linear system at each time step due to the non-diagonal mass matrix.}

At the discrete level, we showed that the implicit method is unconditionally stable \reviewerTwo{in the total energy norm up to the 
Lipschitz continuity of the surface height function} and conserves volume \reviewerTwo{under $\surf>b$ assumption}. 
The semi-implicit scheme, even when stabilized with \reviewerOne{the original} FSSA, is only conditionally stable \reviewerTwo{up to Lipschitz continuity of $\surf$} and does not conserve volume. 
The explicit Euler method conserves volume \reviewerTwo{under $\surf>b$ assumption} but is unstable in the total energy norm. 
\igorr{We theoretically found that the cause of the instabilities is related 
to the normal surface velocity norm. We proposed 
a stabilization term that accurately eliminates these instabilities and showed unconditional stability with respect to $\Delta t$ \reviewerTwo{up to Lipschitz continuity of $\surf$}, while preserving the domain volume \reviewerTwo{under $\surf>b$ assumption}. 
The stabilization term is symmetric, simple to implement, and works without tuning parameters.}
Numerical tests confirmed all the theoretical results. 

The \reviewerTwo{new} stabilized explicit scheme \eqref{eq:forward_euler_stabilized} that we proposed offers a robust and efficient alternative to more expensive implicit solvers. 
This is the scheme we recommend using in practice. 

\reviewerTwo{Understanding what regularization terms are required and computationally practical for well-posedness of the Stokes/free-surface problem 
is an important direction for future work.
} \reviewerOne{Moreover, similar approach towards stability analysis and stabilization design may be possible to apply to 
other free-surface problems such as the Navier-Stokes/free-surface problem \cite{Audusse_FEM_Navier_stokes_Free_Surface}, 
and to other time integrators such as higher-order Runge-Kutta methods.}

\section{Acknowledgements}
Josefin Ahlkrona and Igor Tominec were funded by the Swedish Research Council, grant number 2021-04001.
Lukas Lundgren was funded by the Swedish e-Science Research Centre (SeRC). Igor Tominec 
was also funded by the Slovenian Research and Innovation Agency (ARIS), grant numbers J2-70128 and P2-0162.

\bibliographystyle{plain}
\bibliography{references}


\end{document}